\documentclass[11 pt,reqno]{amsart}
\usepackage{amssymb}
\usepackage{amsmath}
\usepackage{amscd}

\newtheorem {theorem}{Theorem}[section]
\newtheorem {proposition}[theorem]{Proposition}
\newtheorem {corollary}[theorem]{Corollary}
\newtheorem {lemma}[theorem]{Lemma}
\newtheorem {conjecture}[theorem]{Conjecture}

\numberwithin {equation}{section}
\renewcommand {\proof}{{\sc Proof.}\ }

\newcommand {\remark}{{\sc Remark.\ }}

\newcommand {\definition}{{\sc Definition.\ }}

\newcommand {\halmos}{$\blacksquare$}

\newcommand {\IC}{\mathbb{C}}
\newcommand {\IN}{\mathbb{N}}                          
\newcommand {\IP}{\mathbb{P}}
\newcommand {\IQ}{\mathbb{Q}}
\newcommand {\IR}{\mathbb{R}}
\newcommand {\IZ}{\mathbb{Z}}

\newcommand {\A}{\mathcal A}
\newcommand {\C}{\mathcal C}

\newcommand {\F}{\mathcal F}

\renewcommand {\H}{\mathcal H}

\renewcommand {\O}{\mathcal O}
\renewcommand {\P}{\mathcal P}
\newcommand {\R}{\mathcal R}
\newcommand {\U}{\mathcal U}
\newcommand {\V}{\mathcal V}
\newcommand {\X}{\mathcal X}

\renewcommand {\a}{\mathfrak a}
\renewcommand {\b}{\mathfrak b}

\newcommand {\g}{\mathfrak{g}}
\newcommand {\h}{\mathfrak h}
\renewcommand {\ll}{\mathfrak l}
\newcommand {\n}{\mathfrak n}
\newcommand {\p}{\mathfrak p}
\renewcommand {\r}{\mathfrak r}

\renewcommand {\SS}[1]{\mathfrak S_{#1}}

\newcommand {\ad}{\operatorname{ad}}
\newcommand {\Ad}{\operatorname{Ad}}

\newcommand {\End}{\operatorname{End}}
\newcommand {\Hom}{\operatorname{Hom}}
\newcommand {\Ker}{\operatorname{Ker}}

\newcommand {\id}{\operatorname{id}}

\newcommand {\gr}{\operatorname{gr}}

\newcommand {\gl}[1]{\mathfrak{gl}_{#1}}
\renewcommand {\sl}[1]{\mathfrak{sl}_{#1}}

\newcommand {\so}[1]{\mathfrak{so}_{#1}}
\renewcommand {\sp}[1]{\mathfrak{sp}_{#1}}

\newcommand {\Uh}{U\h}
\newcommand {\Ug}{U\g}
\newcommand {\Un}{U\n}
\newcommand {\Ugh}{U\g^{\h}}

\newcommand {\Uhg}{U_{\hbar}\g}

\newcommand {\Uhsl}[1]{U_{\hbar}\sl{#1}}

\newcommand {\Sn}{S\n}
\newcommand {\Sgh}{S\g^{\h}}

\newcommand {\Cg}{\C_{\g}}
\newcommand {\Csln}{\C_{\sl{n}}}

\newcommand {\Csl}[1]{\C_{\sl{#1}}}

\newcommand {\kalpha}{\kappa_{\alpha}}

\newcommand {\calpha}{C_{\alpha}}
\newcommand {\nablac}{\nabla_{C}}
\newcommand {\nablak}{\nabla_{\kappa}}

\newcommand {\reg}{_{\operatorname{reg}}}
\newcommand {\hreg}{\h\reg}
\newcommand {\Pg}{P_{\g}}
\newcommand {\Bg}{B_{\g}}

\newcommand {\ie}{{\it i.e., }}
\newcommand {\eg}{{\it e.g.}, }
\newcommand {\fd}{finite--dimensional }
\newcommand {\fg}{finitely--generated }

\newcommand {\ol}{\overline}
\newcommand {\ul}{\underline}
\newcommand {\wt}{\widetilde}
\newcommand {\wh}{\widehat}

\newcommand {\<}{\langle}
\renewcommand {\>}{\rangle}

\newcommand {\GZ}{Gelfand--Zetlin }
\newcommand {\KR}{Kirillov--Reshetikhin }
\newcommand {\KZ}{Knizhnik--Zamolodchikov  }
\newcommand {\KKZ}{_{\scriptscriptstyle{\operatorname{KZ}}}}
\newcommand {\KT}{Koike--Terada }

\newcommand {\fml}{[\negthinspace[\hbar]\negthinspace]}
\newcommand {\fmll}{[\negthinspace[h]\negthinspace]}

\newcommand {\ICh}{\IC\fml}
\newcommand {\La}{((h))}

\newcommand {\Pus}[1]{\IC((#1^{\frac{1}{\infty}}))}

\newcommand {\Gr}{\operatorname{Gr}}
\newcommand {\red}{\operatorname{red}}

\newcommand {\res}{\operatorname{res}}

\newcommand {\nn}{^{(n)}}
\newcommand {\Th}{\Theta_{\hbar}}
\newcommand {\mmu}[2]{\mu^{(#1)}_{#2}}

\newcommand {\half}[1]{\frac{#1}{2}}

\newcommand {\sqnm}[1]{\<#1,#1\>}
\newcommand {\cor}[1]{\alpha_{#1}^{\vee}}

\newcommand {\Gal}{\operatorname{Gal}}

\input{diagrams}

\begin{document}

\title
[Casimir Operators and Generalised Braid Groups]
{Casimir Operators and\\
Monodromy Representations of\\
Generalised Braid Groups}
\author[J. J. Millson]{John J. Millson$^*$}
\address{
Mathematics Department \\
University of Maryland \\
College Park, MD 20742--4015}
\email{jjm@math.umd.edu}
\thanks
{$^*$ work partially supported by NSF grants
DMS-98-03520 and 01-04006.}
\author[V. Toledano Laredo]{Valerio Toledano Laredo$^{**}$}
\address{
Institut de Mathematiques de Jussieu \\
Universite Pierre et Marie Curie     \\
UMR 7586, Case 191                   \\
16 rue Clisson\\
F--75013 Paris}
\email{toledano@math.jussieu.fr}
\thanks
{$^{**}$  work partially supported by an MSRI
postdoctoral fellowship for the year 2000--2001}
\date{September 2004}
\begin{abstract}
Let $\g$ be a complex, simple Lie algebra with Cartan
subalgebra $\h$ and Weyl group $W$. We construct a
one--parameter family of flat connections $\nablak$
on $\h$ with values in any finite--dimensional $\g
$--module $V$ and simple poles on the root hyperplanes.
The corresponding monodromy representation of the braid
group $\Bg$ of type $\g$ is a deformation of the action
of (a finite extension of) $W$ on $V$. The residues of
$\nablak$ are the Casimirs $\kalpha$ of the subalgebras
$\sl{2}^{\alpha}\subset\g$ corresponding to the roots of
$\g$. The irreducibility of a subspace $U\subseteq V$ under
the $\kalpha$ implies that, for generic values of the
parameter, the braid group $\Bg$ acts irreducibly on
$U$. Answering a question of Knutson and Procesi, we
show that these Casimirs act irreducibly on the
weight spaces of all simple $\g$--modules if $\g
=\sl{3}$ but that this is not the case if $\g\ncong
\sl{2},\sl{3}$. We use this to disprove a conjecture
of Kwon and Lusztig stating the irreducibility of
quantum Weyl group actions of Artin's braid group
$B_{n}$ on the zero weight spaces of all simple
$\Uhsl{n}$--modules for $n\geq 4$. Finally, we study
the irreducibility of the action of the Casimirs on
the zero weight spaces of self--dual $\g$--modules
and obtain complete classification results for $\g
=\sl{n}$ and $\g_{2}$.
\end{abstract}
\maketitle

\setcounter{tocdepth}{1}
\tableofcontents

\section{Introduction}\label{se:intro}

It has been known since the seminal work of Knizhnik and
Zamolodchikov how to construct representations of Artin's
braid groups $B_{n}$ by using the representation theory
of a given complex, semi--simple Lie algebra $\g$ \cite{KZ}.
Realising $B_{n}$ as the fundamental group of the quotient
of the configuration space
$$X_{n}=
\{(z_{1},\ldots,z_{n})\in\IC^{n}|
\medspace z_{i}\neq z_{j},\medspace 1\leq i<j\leq n\}$$
by the natural action of the symmetric group $\SS{n}$,
one obtains these representations as the monodromy of
the \KZ connection
$$\nabla\KKZ=d-h\!\!\!\!\sum_{1\leq i<j\leq n}\!\!\!\!
\frac{d(z_{i}-z_{j})}{z_{i}-z_{j}}\cdot\Omega_{ij}$$
with values in the $n$--fold tensor product $V^{\otimes n}$
of a finite--dimensional $\g$--module $V$. Here, the one--form
$\nabla\KKZ$ is regarded as an $\SS{n}$--equivariant flat
connection on the topologically trivial vector bundle over
$X_{n}$ with fibre $V^{\otimes n}$ and then pushed down to
$X_{n}/\SS{n}$. Its coefficients $\Omega_{ij}\in\End(V^{\otimes
n})$ are given by
$$\Omega_{ij}=\sum_{a=1}^{\dim\g}\pi_{i}(X_{a})\pi_{j}(X^{a})$$
where $\{X_{a}\},\{X^{a}\}$ are dual basis of $\g$ with respect
to the Killing form and $\pi_{k}(\cdot)$ denotes the action on
the $k$th tensor factor of $V^{\otimes n}$. Finally, the complex
number $h$ may be regarded as a deformation parameter which, upon
being set to 0, gives a monodromy representation of $B_{n}$
factoring through the natural action of the symmetric group
on $V^{\otimes n}$.\\

Aside from their intrinsic interest, these representations
appear naturally in a number of different contexts. They
define for example the commutativity and associativity
constraints in the tensor category of highest weight
representations of the affine Kac--Moody algebra $\wh{\g}$
\cite{KL,Wa,TL1} and, by the Kohno--Drinfeld theorem, on
the finite--dimensional representations of the quantum
group $\Uhg$ \cite{Dr3,Dr4,Ko1}. As such, they define
invariants of knots and links and, for suitable rational
values of $h$, of three--manifolds \cite{Tu}.\\

The purpose of the present paper is to use the representation
theory of $\g$ in a similar vein to construct monodromy
representations of a different braid group, namely the
generalised braid group $\Bg$ of type $\g$. The latter
may be defined as the fundamental group of the quotient
$\hreg/W$ of the set $\hreg$ of regular elements in a
Cartan subalgebra $\h$ of $\g$ by the action of the
corresponding Weyl group $W$. Like Artin's braid groups,
$\Bg$ is presented on generators $S_1,\ldots,S_n$
labelled by a choice $\alpha_{1},\ldots,\alpha_{n}$ of
simple roots of $\g$ with relations
$$S_iS_j\cdots=S_jS_i\cdots$$
for any $i\neq j$, where the number of factors on each
side is equal to the order of the product $s_{i}s_{j}$
of the orthogonal reflections corresponding to $\alpha
_{i}$ and $\alpha_{j}$ in $W$ \cite{Br}.\\

To state our first main result, let $R=\{\alpha\}\subset\h^{*}$
be the set of roots of $\g$ relative to $\h$ so that $\hreg
=\h\setminus\bigcup_{\alpha\in R}\Ker(\alpha)$. For each
$\alpha\in R$, let $\sl{2}^{\alpha}=\<e_{\alpha},f_{\alpha},
h_{\alpha}\>\subseteq\g$ be the corresponding $\sl{2}(\IC)
$--subalgebra of $\g$ and
$$\kalpha=
\frac{\<\alpha,\alpha\>}{2}
(e_{\alpha}f_{\alpha}+f_{\alpha}e_{\alpha})$$
the truncated Casimir operator of $\sl{2}^{\alpha}$ where
$\<\cdot,\cdot\>$ is a fixed multiple of the Killing form
of $\g$. Let $V$ be a $\g$--module, then we prove in section
\ref{se:flat} the following\footnote{Theorem \ref{th:i nablak}
was independently discovered by De Concini around 1995
(unpublished). A variant of the connection $\nablak$ also
appears in the recent paper \cite{FMTV}. Unlike $\nablak$
however, the connection introduced in \cite{FMTV} is not
$W$--equivariant and therefore only defines representations
of the {\it pure} braid group $\Pg$ instead of the full braid
group $\Bg$}

\begin{theorem}\label{th:i nablak}
The one--form
$$\nablak=
d-h\sum_{\alpha\in R}\frac{d\alpha}{\alpha}\cdot\kalpha$$
defines, for any $h\in\IC$, a flat connection on the
trivial vector bundle over $\hreg$ with fibre $V$,
which is reducible with respect to the weight space
decomposition of $V$.
\end{theorem}

As a consequence, each weight space of $V$ carries a canonical
one--parameter family of monodromy representations of the {\it
pure} braid group $\Pg=\pi_{1}(\hreg)$. This action extends to
one of the full braid group $\Bg$ on the direct sum of weight
spaces corresponding to a given Weyl group orbit, and in
particular on the zero weight space of $V$, by pushing $\nablak
$ down to the quotient space $\hreg/W$. Since the Weyl group
itself does not act on $V$, this requires choosing an action
of $\Bg$ on $V$ which permutes the weight spaces compatibly
with the projection $\Bg\rightarrow W$. This may for example
be achieved by taking the simply--connected complex Lie group
$G$ corresponding to $\g$ and mapping $\Bg$ to one of the
Tits extension $\wt W$ of $W$, a class of subgroups of the
normaliser in $G$ of the torus $T$ corresponding to $\h$
which are extensions of $W$ by the sign group $\IZ_{2}^{n}$,
where $n=\dim(\h)$ \cite {Ti}. The choice of a specific Tits
extension is somewhat immaterial since any two are conjugate
by an element of $T$ and the corresponding representations
of $\Bg$ are therefore equivalent.\\

The rest of the paper is devoted to the study of the
irreducibility of our monodromy representations. Define
a subspace $U\subseteq V$ invariant under the monodromy
action of $\Bg$ to be {\it generically irreducible} if it is irreducible
for all values of $h$ lying outside the zero set of some holomorphic
function. In section \ref{se:gen irred}, we prove the following 

\begin{theorem}\label{th:i gen irred}
A subspace $U\subseteq V$ is generically irreducible under
the braid group $\Bg$ (resp. the pure braid group $\Pg$)
if, and only if it is irreducibly acted upon by the Casimirs
$\kalpha$ and $\wt W$ (resp. the $\kalpha$ and $\wt{W}\cap T$).
\end{theorem}

This naturally prompts the question, originally asked us by
C. Procesi and A. Knutson, of whether the {\it Casimir algebra}
\ie the algebra
$$\Cg=\<\kalpha\>_{\alpha\in R}\vee\h\subset\Ug$$
generated by the Casimirs $\kalpha$ and $\h$ inside the
enveloping algebra of $\g$ acts irreducibly on the weight
spaces of any simple $\g$--module or, stronger still, whether
 it is equal to the algebra $\Ugh$ of $\h$--invariants in
$U\g$.\\

The answer to both questions is clearly positive for $\g=
\sl{2}$ and we show in section \ref{se:Casimir} that it
this almost so for $\g=\sl{3}$. More precisely,

\begin{theorem}\label{th:i Casimir sl3}
If $\g=\sl{3}$, $\Cg$ is a proper subalgebra of $\Ugh$, but
the latter is generated by $\Cg$ and the centre $Z(\Ug)$ of
$\Ug$. In particular, $\Cg$ acts irreducibly on the weight
spaces of any simple $\g$--module.
\end{theorem}

As a consequence, all monodromy representations of $P_{3}$
on weight spaces of simple $\sl {3}$--modules, and of $B_{3}$
on their zero weight spaces are generically irreducible, a
fact which refines a result proved by Kwon \cite{Kw} in the
context of quantum Weyl groups, and to which we shall return
below. For $\g\ncong\sl{2},\sl{3}$, the situation is radically
different and we prove

\begin{theorem}\label{th:i Casimir notsl3}
If $\g\ncong\sl{2},\sl{3}$, there exists a simple $\g$--module
$V$ the zero weight space of which is reducible under the
joint action of $\Cg$ and of $W$. In particular, $\Cg$ and
$Z(U\g)$ do not generate $\Ugh$.
\end{theorem}

For $\g\ncong\sl{n}$, our $V$ is in fact the kernel of the
commutator map $[\cdot,\cdot]:\g\wedge\g\rightarrow\g$.
For $\g\cong\sl{n}$, $\Ker([\cdot,\cdot])$ is reducible
and the construction of a suitable $V$ relies on the
following general reducibility criterion, valid for any
$\g$. Let $V$ be a simple $\g$--module with zero weight
space $V[0]\neq\{0\}$. If $V$ is self--dual, it is
acted on by a linear involution $\Theta_{V}$ such that,
for any $X\in\g$, $\Theta_{V}X\Theta_{V}^{-1}=\Theta(X)$
where $\Theta$ is the Chevalley involution of $\g$ relative
to a given choice of simple root vectors. Since $\Theta$
acts as $-1$ on $\h$ and fixes the Casimirs $\kalpha$ and
$\wt W$, $\Theta_{V}$ leaves $V[0]$ invariant and commutes 
with $\Cg$ and $\wt W$. $V[0]$ is therefore reducible under
$\Cg$ and $\wt{W}$ whenever $\Theta_{V}$ does not act as a
scalar on it.\\

To prove that this is the case for some $V$, we note further
that if $\r\subset\g$ is a reductive subalgebra normalised
by $\Theta$ and $V$ is such that its restriction to $\r$
contains a zero--weight vector $u$ lying in a simple $\r
$--summand $U$ which is not self--dual, then $\Theta_{V}u$
cannot be proportional to $u$ since $U\cap\Theta_{V} U=\{0\}$.
To summarise, our initial problem reduces to finding
simple, self--dual $\g$--modules $V$ whose restriction
to some reductive subalgebra $\r\subset\g$ contains
non--self dual summands intersecting $V[0]$ non--trivially.
For $\g=\sl{n}$, we construct such $V$'s by using the
\GZ branching rules for the inclusion $\gl{n-1}\subset
\gl{n}$ \cite{GZ}.\\

In section \ref{se:Kwon}, we use our results to disprove
a conjecture of Kwon and Lusztig on quantum Weyl group
actions of the braid group $B_{n}$ \cite{Kw}. To state
it, recall that the Drinfeld--Jimbo quantum group $\Uhg$
corresponding to $\g$ defines, on any of its integrable
representations $\V$, an action of the braid group $\Bg$
called the {\it quantum Weyl group action}, which is a
deformation of the action of $\wt W$ on the $\g$--module
$V=\V/\hbar\V$ \cite{Lu,KR,So}. In \cite{Kw}, Kwon considered
the case of $\g=\sl{n}$ and gave a necessary condition
for the zero weight space of $\V$ to be irreducible under
$\Bg=B_{n}$. He showed in particular that the zero weight
spaces of all $U_{\hbar}\sl{3}$--modules are irreducible
under $B_{3}$. Based on these findings he and Lusztig
conjectured that this should hold for all $B_{n}$, $n
\geq 4$. 

\begin{theorem}\label{th:i Kwon}
The Kwon--Lusztig conjecture is false for any simple,
complex Lie algebra $\g\ncong\sl{2},\sl{3}$.
\end{theorem}

Our disproof is based on the simple observation that the
{\it quantum} Chevalley involution $\Th$ of $\Uhg$ acts
on any self--dual $\Uhg$--module $\V$ and that its
restriction to the zero weight space $\V[0]$ centralises
the action of $\Bg$.  We then remark that $\Th$ acts as
a scalar on $\V[0]$ iff the classical Chevalley involution
acts as a scalar on the zero weight space of the
$\g$--module $V=\V/\hbar\V$ and rely on the results
of section \ref{se:Casimir}.\\

In section \ref{se:irreps} we show that, despite the
reducibility results of \S \ref{se:Casimir}, the connection
$\nablak$ yields irreducible monodromy representations
of $\Bg$ of arbitrarily large dimensions. For $\g\cong\sl{n}$, 
we show in fact that the weight spaces of all Cartan powers
of the adjoint representation are irreducible under the the
Casimirs $\kalpha$.\\

In section \ref{se:selfdual V[0]} we show that, when
$\g$ is isomorphic to $\sl{n}$, $n\geq 4$, or $\g_{2}$,
the zero weight space of most {\it self--dual}, simple
$\g$--modules is reducible under the Casimir algebra
$\Cg$ of $\g$, thus strengthening the results of section
\ref{se:Casimir}. More precisely, let $V$ be a non--trivial,
simple, self--dual $\g$--module with zero weight space
$V[0]\neq\{0\}$. Then, we obtain the following classification
results

\begin{theorem}\label{th:intro sl}
If $\g\cong\sl{n}$, $n\geq 4$, $V[0]$ is irreducible under
$\Cg$ if, and only if its highest weight $\lambda\in\IZ^n$
is of one  of the following forms
\begin{enumerate}
\item $\lambda=(p,0,\ldots,0,-p)$, $p\in\IN$.
\item $\lambda=(\underbrace{1,\ldots,1}_{k},0,\ldots,0,
\underbrace{-1,\ldots,-1}_{k})$, $0\leq k\leq n/2$.
\item $\lambda=(p,p,-p,-p)$, $p\in\IN$.
\end{enumerate}
\end{theorem}

\begin{theorem}\label{th:intro g2}
If $\g\cong\g_{2}$, $V[0]$ is irreducible under $\Cg$ if,
and only if $V$ is fundamental representation or its
second Cartan power.
\end{theorem}

Partial classification results pointing to the same phenomenon
are obtained in \cite{HMTL} for $\g=\so{2n},\so{2n+1},\sp{2n}$.
The proofs of theorems \ref{th:intro sl}--\ref{th:intro g2} rely on
the use of the Chevalley involution $\Theta$ outlined above
and branching to the subalgebras $\gl{k}\subset\sl{n}$ and
$\sl{3}\subset\g_{2}$ respectively. They show in fact that $
V[0]$ is irreducible under $\Cg$ if, and only if $\Theta$ acts
as a scalar on it. It seems natural to conjecture that this should
be so for any $\g$.\\

It is interesting to note how the reducibility results
of section \ref{se:selfdual V[0]} contrast with the
following theorem of Etingof, which is reproduced with
his kind permission in section \ref{se:etingof}. Let
$\beta\in\sum_{i=1}^{n}\IN\cdot\alpha_{i}$ be a fixed,
non--negative linear combination of simple roots and,
for $\mu\in\h^*$, let $M_{\mu}[\mu-\beta]$ be the subspace
of weight $\mu-\beta$ of the Verma module of highest
weight $\mu$.

\begin{theorem}[Etingof]
There exists a Zariski open set $\O_{\beta}\subset
\h^{*}$ such that, for any $\mu\in\O_{\beta}$, $M_
{\mu}[\mu-\beta]$ is irreducible under the Casimir
algebra $\Cg$.
\end{theorem}

The above theorem has as corollary the following
interesting result, which is also proved in \S \ref
{se:etingof}

\begin{theorem}[Etingof]
The centraliser of the Casimir algebra $\Cg$ in $\Ug$
is generated by $\h$ and the centre of $\Ug$.
\end{theorem}

{\bf Acknowledgements.} We wish to heartily thank Allen
Knutson whose observations on reading \cite{KM} initiated
this project. We are also grateful to R. Buchweitz, C. De
Concini, P. Etingof, M. Kashiwara, F. Knop, B. Kostant, C.
Laskowski, A. Okounkov, C. Procesi, R. Rouquier, S. Yuzvinksy
and A. Wassermann for a number of useful discussions
and the referees for pointing out a number of inaccuracies
and improvements to the exposition in
an earlier version of this paper. The research for this
project was partly carried out while the first author
visited the Institut de Math\'ematiques de Jussieu in
June 2000 and 2001 and while the second author
was a post--doctoral fellow at MSRI during the academic
year 2000--2001. We are grateful to both institutions
for their financial support and pleasant working conditions.

\section{Flat connections on $\h\reg$}\label{se:flat}

\subsection{The flat connection $\nablak$}\label{ss:nablak}

Let $\g$ be a complex, semi--simple Lie algebra with Cartan
subalgebra $\h$ and root system $R=\{\alpha\}\subset\h^{*}$.
Let
$$\hreg=\h\setminus\bigcup_{\alpha\in R}\Ker(\alpha)$$
be the set of regular elements in $\h$ and $V$ a \fd
$\g$--module. We shall presently define a flat connection
on the trivial vector bundle $\hreg\times V$ over $\hreg
$\footnote{all vector bundles and connections considered
in this section are holomorphic.}. We need for this purpose
the following flatness criterion due to Kohno \cite{Ko2}. Let
$B$ be a complex, finite--dimensional vector space and
$\A=\{H_{i}\}_{i\in I}$ a finite collection of hyperplanes in
$B$ determined by the linear forms $\phi_{i}\in B^{*}$, $i\in
I$.

\begin{lemma}\label{le:kohno}
Let $V$ be a finite--dimensional vector space and $\{r_{i}\}
\subset\End(V)$ a family indexed by $I$. Then,
\begin{equation}\label{eq:kohno}
\nabla=d-\sum_{i\in I}\frac{d\phi_{i}}{\phi_{i}}\cdot r_{i}
\end{equation}
defines a flat connection on $(B\setminus\A)\times V$
iff, for any subset $J\subseteq I$ maximal for the
property that $\bigcap_{j\in J}H_{j}$ is of codimension
2, the following relations hold for any $j\in J$
\begin{equation}\label{eq:relations}
[r_{j},\sum_{j'\in J}r_{j'}]=0
\end{equation}
\end{lemma}

\remark Since the relations \eqref{eq:relations} are homogeneous,
a solution $\{r_{i}\}_{i\in I}$ of \eqref {eq:relations} defines
in fact a {\it one--parameter} family of representations
$$\rho_{h}:\pi_{1}(B\setminus\A)\longrightarrow GL(V)$$
parametrised by $h\in\IC$ where $\rho_{h}$ is the monodromy
of the connection \eqref{eq:kohno} with $r_{i}$ replaced by
$h\cdot r_{i}$.\\

For any $\alpha\in R$, choose root vectors $e_{\alpha}\in
\g_{\alpha},f_{\alpha}\in\g_{-\alpha}$ such that
$[e_{\alpha},f_{\alpha}]=h_{\alpha}=\alpha^{\vee}$ and let
\begin{equation}\label{eq:truncated}
\kappa_{\alpha}=
\half{\sqnm{\alpha}}
(e_{\alpha}f_{\alpha}+f_{\alpha}e_{\alpha})
\in U\g
\end{equation}
be the truncated Casimir operator of the three--dimensional
subalgebra $\sl{2}^{\alpha}\subset\g$ spanned by $e_{\alpha}
,h_{\alpha},f_{\alpha}$ relative to the restriction to $\sl
{2}^{\alpha}$ of a fixed multiple $\<\cdot,\cdot\>$ of the
Killing form of $\g$. Note that $\kappa_{\alpha}$ does not
depend upon the particular choice of $e_{\alpha}$ and $f_
{\alpha}$ and that $\kappa_{-\alpha}=\kappa_{\alpha}$. Let
$R^{+}\subset R$ be the set of positive roots corresponding
to a choice of simple roots $\alpha_{1},\ldots,\alpha_{n}$
of $\g$. 

\begin{theorem}\label{th:casimir flat}
The one--form
\begin{equation}\label{eq:Casimir connection}
\nablak=
d-h\sum_{\alpha\in R^{+}}
\frac{d\alpha}{\alpha}\cdot\kalpha=
d-\frac{h}{2}\sum_{\alpha\in R}
\frac{d\alpha}{\alpha}\cdot\kalpha
\end{equation}
defines, for any $h\in\IC$, a flat connection on
$\hreg\times V$ which is reducible with respect
to the weight space decomposition of $V$.
\end{theorem}
\proof By lemma \ref{le:kohno}, we must show that
for any rank 2 root subsystem $R_{0}\subseteq R$
determined by the intersection of $R$ with a 2--dimensional
subspace in $\h^*$, the following holds for any
$\alpha\in R_{0}^{+}=R_{0}\cap R^{+}$
\begin{equation}\label{eq:rank 2}
[\kappa_{\alpha},\sum_{\beta\in R_{0}^{+}}\kappa_{\beta}]=0
\end{equation}
This may be proved by an explicit computation by
considering in turn the cases where $R_{0}$ is of
type $A_{1}\times A_{1}$, $A_{2}$, $B_{2}$ or
$G_{2}$ but is more easily settled by the following
elegant observation of A. Knutson \cite{Kn}. Let
$\g_{0}\subseteq\g$ be the semi--simple Lie algebra
with root system $R_{0}$, $\h_{0}\subset\h$ its Cartan
subalgebra and $C_{0}\in Z(\Ug_{0})$ its Casimir
operator. Then, $\sum_{\beta\in R_{0}^{+}}\kappa_{\beta}
-C_{0}$ lies in $U\h_{0}$ so that \eqref{eq:rank 2}
holds since $\kappa_{\alpha}$ commutes with $\h_{0}$.
The reducibility of $\nablak$ with respect to the
$\h$--action on $V$ is an immediate consequence of
the fact that the operators $\kalpha$ are of weight
zero \halmos\\

\remark Altough $V$ admits a hermitian inner product with
respect to which the Casimirs $\kalpha$ are self--adjoint,
it is easy to check that the connection $\nablak$ is not
unitary with respect to the corresponding constant inner
product on $\hreg\times V$. However, the fact that the
connection $\nablak$ for $\g=\sl{n}$ coincides with the
(genus 0) \KZ connection on $n$ points for $\g'=\sl{k}$
via Howe duality \cite[thm. 3.5]{TL2}, and that the latter
is conjectured to be unitary on the subbundle of conformal
blocks for suitable rational values of $h$ \cite{Ga}
\footnote{this is now a theorem, at least for $\g=\sl{2}$,
see \cite{Ra}}, suggests that the connection $\nablak$
ought to be unitary for any $\g$. It is an interesting
open problem to determine whether this is so.\\

Let $W$ be the Weyl group of $\g$. Fix a basepoint $t_0
\in\hreg$, let $[t_0]$ be its image in $\hreg/W$ and let $
\Pg=\pi_{1}(\hreg;t_0)$, $\Bg=\pi_{1}(\hreg/W;[t_0])$ be the
generalised pure and full braid groups of type $\g$. We
adhere here to the convention that for a pointed topological
space $(X,x_0)$, the composition $q\cdot p$ of two paths
$p,q\in\pi_1(X;x_0)$ is given by $p$ {\it followed} by $q$
so that the holonomy of a flat vector bundle $(\V,\nabla)$
at $x_0$ is a group homomorphism $\pi(X;x_0)\rightarrow
GL(\V_{x_0})$. The fibration $\hreg\rightarrow\hreg/W$
gives rise to the exact sequence
$$1\longrightarrow\Pg\longrightarrow\Bg
\longrightarrow W\longrightarrow 1$$
where the rightmost arrow is obtained by associating to
$p\in\Bg$ the unique $w\in W$ such that $w^{-1}t_0=
\wt{p}(1)$ where $\wt{p}$ is the unique lift of $p$ to a
path in $\hreg$ such that $\wt{p}(0)=t_0$.\\

By theorem \ref{th:casimir flat}, the monodromy of $\nablak$
yields a one--parameter family of representations of $\Pg$
on $V$ preserving its weight space decomposition. We wish
to extend this action to one of $\Bg$, by pushing $\nablak$
down to a flat connection on the quotient $\hreg/W$. Since
$W$ does not act on $V$, this requires choosing an action
of $\Bg$ on $V$. Let for this purpose $G$ be the complex,
connected and simply--connected Lie group with Lie algebra
$\g$, $T$ its torus with Lie algebra $\h$ and $N(T)\subset
G$ the normaliser of $T$ so that $W\cong N(T)/T$. We regard
$\Bg$ as acting on $V$ by choosing a homomorphism $\sigma:
\Bg\rightarrow N(T)$ compatible with
\begin{equation}\label{eq:extension}
\begin{diagram}[height=2.5em]
\Bg&\rTo^{\sigma}&N(T) \\
   &\rdTo        &\dTo \\
   &             &W
\end{diagram}
\end{equation}
Such $\sigma$'s abound and we describe in \S \ref{ss:Tits}
a class of them which we call Tits extensions \cite{Ti}.
Let $\wt{\hreg}\xrightarrow {p}\hreg$ be the universal
cover of $\hreg$ and $\hreg/W$.

\begin{proposition}\label{th:existence}
The one--form $p^{*}\nablak$ defines a $\Bg$--equivariant
flat connection on $p^{*}(\hreg\times V)=\wt\hreg\times V$.
It therefore descends to a flat connection on the vector bundle
$$\begin{diagram}[height=2em,width=2em]
V&\rTo&\wt{\hreg}\times_{\Bg}V \\
 &    &\dTo                    \\
 &    &\hreg/W
\end{diagram}$$
which is reducible with respect to the weight space
decomposition of $V$.
\end{proposition}
\proof $\Bg$ acts on $\Omega^{\bullet}(\wt\hreg,V)
=\Omega^{\bullet}(\wt\hreg)\otimes V$ by $\gamma
\rightarrow(\gamma^{-1})^{*}\otimes\sigma(\gamma)$.
Thus, if $\gamma\in\Bg$ projects onto $w\in W$, we
get using $p\cdot\gamma^{-1}=w^{-1}\cdot p$,
$$\gamma\medspace p^{*}\nablak\medspace\gamma^{-1}=
d-
\half{h}\sum_{\alpha\in R}dp^{*}w\alpha/p^{*}w\alpha
\otimes\sigma(\gamma)\kappa_{\alpha}\sigma(\gamma)^{-1}$$
Since $\kalpha$ is independent of the choice of the root
vectors $e_{\alpha},f_{\alpha}$ in \eqref{eq:truncated},
$\Ad(\sigma(\gamma))\kappa_{\alpha}=\kappa_{w\alpha}$
and the above is equal to $p^{*}\nablak$ as claimed. $p^{*}
\nablak$ is flat and commutes with the fibrewise action of
$\h$ by theorem \ref{th:casimir flat} \halmos\\

Thus, for any homomorphism $\sigma:\Bg\rightarrow N(T)$
compatible with \eqref{eq:extension}, proposition
\ref{th:existence} yields a one--parameter family
of monodromy representations
$$\rho^{\sigma}_{h}:\Bg\longrightarrow GL(V)$$
which permutes the weight spaces of $V$ compatibly
with the action of $W$ on $\h^*$. By standard ODE
theory, $\rho^{\sigma}_{h}$ depends analytically on the
complex parameter $h$ and, when $h=0$, is equal to
the action of $\Bg$ on $V$ given by $\sigma$. We
record for later use the following elementary

\begin{proposition}\label{pr:recipe}
Let $\gamma\in\Bg=\pi_{1}(\hreg/W;[t_0])$ and let
$\wt\gamma:[0,1]\rightarrow\hreg$ be a lift of $\gamma$.
Then, 
\begin{equation}\label{eq:recipe}
\rho^{\sigma}_{h}(\gamma)=\sigma(\gamma)\P(\wt\gamma)
\end{equation}
where $\P(\wt\gamma)\in GL(V)$ is the parallel
transport along $\wt\gamma$ for the connection
$\nablak$ on $\hreg\times V$.
\end{proposition}
\proof
Let $\wt{\wt\gamma}:[0,1]\rightarrow\wt{\hreg}$
be a lift of $\gamma$ and $\wt{\gamma}$ so that
$\wt{\wt\gamma}(1)=\gamma^{-1}\wt{\wt\gamma}(0)$.
Then, since the connection on $p^{*}(\hreg\times
V)$ is the pull--back of $\nablak$, and that
on $\left(p^{*}(\hreg\times V)\right)/\Bg$ the
quotient of $p^{*}\nablak$, we find that
$$\rho^{\sigma}_{h}(\gamma)=
\P(\gamma)=
\sigma(\gamma)\P(\wt{\wt\gamma})=
\sigma(\gamma)\P(\wt\gamma)$$
\halmos\\

\remark By \eqref{eq:recipe}, the representation $\rho^
{\sigma}_{h}$ depends on the choice of the homomorphism
$\sigma$ satisfying \eqref{eq:extension}. However, the
representations corresponding to Tits extensions are all
equivalent since these are conjugate by elements of $T$ (see
\S \ref{ss:Tits}). Note also that the restriction of $\rho_{h}^
{\sigma}$ to the zero weight space of $V$ does not depend
on $\sigma$ since $W\cong N(T)/T$ acts canonically on it.

\remark Note that, by \eqref{eq:recipe}, the restriction
of $\rho^{\sigma}_{h}$ to the pure braid group $\Pg$
does not coincide with the monodromy of the connection
$\nablak$. Rather, it differs from it by the $T$--valued
character given by the restriction of $\sigma$ to $\Pg$.\\

\remark By Brieskorn's theorem, $\Bg$ is presented on
generators $S_{1},\ldots,S_{n}$ labelled by the simple
simple reflections $s_{1},\ldots,s_{n}\in W$ with relations
\begin{equation}\label{eq:braid relations}
\underbrace{S_{i}S_{j}\cdots}_{m_{ij}}=
\underbrace{S_{j}S_{i}\cdots}_{m_{ij}}
\end{equation}
for any $1\leq i<j\leq n$ where the number $m_{ij}$ of
factors on each side is equal to the order of $s_{i}s_{j}$
in $W$ \cite{Br}. Each $S_{i}$ may be obtained as a small
loop in $\hreg/W$ around the reflecting hyperplane $\Ker
(\alpha_{i})$ of $s_{i}$.

\subsection{Variants of $\nablak$}

If $p_{\alpha}\in\Uh$, $\alpha\in R$, is a collection of
polynomials in $\h$, the connection
$$d-\half{h}\sum_{\alpha\in R}
\frac{d\alpha}{\alpha}\cdot(\kalpha+p_{\alpha})$$

is flat by theorem \ref{th:casimir flat} since $[\kalpha,
p_{\beta}]=[p_{\alpha},p_{\beta}]=0$ for any $\alpha,\beta
\in R$. It is moreover $W$--equivariant if, in addition,
$w p_{\alpha}=p_{w\alpha}$ for any $w\in W$. The corresponding
monodromy representation of $\Pg$ is equal to that of the
connection $\nablak$ tensored with the character $\chi:
\Pg\rightarrow T$ given by the monodromy of the abelian
connection
$$d-\half{h}\sum_{\alpha\in R}\frac{d\alpha}{\alpha}\cdot p_{\alpha}$$
and therefore does not significantly differ from the
monodromy of $\nablak$. A possible choice is to set
$p_{\alpha}=\<\alpha,\alpha\>/2\cdot\medspace h_{\alpha}^{2}$
which yields the connection
$$\nablac=d-h\sum_{\alpha\in R^+}\frac{d\alpha}{\alpha}\cdot\calpha$$
where $C_{\alpha}\in U\sl{2}^{\alpha}$ is the full Casimir
operator of $\sl{2}^{\alpha}$.

\subsection{The holonomy Lie algebra $\a(\A)$}

Kohno's lemma \ref{le:kohno} gives a description of the
holonomy Lie algebra $\a(\A)$ of a general hyperplane
arrangement $\A=\{\H_{i}\}_{i\in I}$ as the quotient of
the free Lie algebra on generators $\{r_{i}\}_{i\in I}$
by the relations \eqref{eq:relations}. When $\A=\A_{\g}
=\{\Ker(\alpha)\}_{\alpha\in R}$ is the arrangement of
root hyperplanes of $\g$, theorem \ref{th:casimir flat}
is equivalent to the fact that the assignment $r_{\alpha}
\rightarrow\kalpha$ extends to an algebra homomorphism
$$\phi:U\a(\A_{\g})\longrightarrow U\g$$
of the universal enveloping algebra of $\a(\A_{\g})$ to that
of $\g$ satisfying
$$\phi(U\a(\A_{\g})_{m})\subset U\g_{2m}$$
for any $m\in\IN$, where the superscript denotes the degree
corresponding to the natural filtrations on both algebras.
We simply note here the following

\begin{proposition}\label{pr:not inj}
If one of the simple factors of $\g$ is not isomorphic to
$\sl{2}$, the map $\phi:U\a(\A_{\g})\longrightarrow U\g$
is not injective.
\end{proposition}
\proof The following argument was pointed out to us by R.
Buchweitz. It suffices to show that, if $\g\ncong\sl{2}$,
$\a(\A_{\g})$ contains a free Lie algebra on at least two
generators, for
then $U\a(\A_{\g})$ has exponential growth with respect to
its filtration, whereas $U\g$, being isomorphic to $S\g$,
only grows polynomially. 
Let $A_{1}\times A_{1}\ncong R_{0}\subseteq R$ be a rank
two root subsystem with positive roots $\beta_{1},\ldots,
\beta_{p}$, $p\geq 3$. Let $\F_{p-1}$ be the free Lie
algebra on generators $x_{1},\ldots,x_{p-1}$ and consider
the maps
$$\F_{p-1}\xrightarrow{i}\a(\A_{\g})\xrightarrow{\pi}\F_{p-1}$$
given by $i(x_{j})=r_{\beta_{j}}$, $j=1\ldots p-1$ and
$$\pi(r_{\alpha})=
\left\{\begin{array}{cl}
0&\text{if $\alpha\notin R_{0}$}\\
x_{j}&\text{if $\alpha=\beta_{j}$, with $1\leq j\leq p-1$}\\
-\sum_{j=1}^{p-1}x_{j}&\text{if $\alpha=\beta_{p}$}
\end{array}\right.$$
It is easy to see that $\pi$ is well--defined, so that $\pi
\circ i=\id$ and $i$ gives an embedding of $\F_{p-1}$ into
$\a(\A_{\g})$ \halmos\\

\remark It seems an interesting problem to find a generating
set of relations for the kernel of the map $\phi$ above. One
such relation may be obtained for any rank two root subsystem
$R_{2}\subset R$ such that the intersection of its $\IZ$--span
with $R$ is equal to $R_{2}$ but the intersection of its $\IR
$--span with $R$ strictly contains $R_{2}$. This is the case
for the root system of type $A_{2}$ given by the long roots
in the root system of type $\g_{2}$ or for the root system
of type $A_{1}\times A_{1}$ generated by any pair of long
roots in the root system of type $C_{n}$. One then has $[
\kalpha,\sum_{\beta\in R_{2}}\kappa_{\beta}]=0$ in $\Ug$,
but $[r_{\alpha},\sum_{\beta\in R_{2}}r_{\beta}]\neq 0$
in $\a(\A_{\g})$.

\subsection{Triviality of $\wt{\hreg}\times_{\Pg}V$}

We show below that the pull--back to $\hreg$ of the vector
bundle $\wt{\hreg}\times_{\Bg}V$ constructed in proposition
\ref{th:existence}, namely
$$\V=\wt{\hreg}\times_{\Pg}V$$
is trivial\footnote{the second author is grateful to R. Rouquier
for a long walk in the Berkeley hills, during which we took turns
in convincing each other that the bundle was trivial, then
non--trivial, then trivial again, until sheer exhaustion and
the late hour of the night suspended, but did not resolve,
the argument.}. Our trivialisation of $\V$ yields an action
of $W$ on $\hreg\times V$ by a one--cocycle, that is by
$w\medspace(t,v)=(wt,A(w,t)v)$ where $A(w,t)\in GL(V)$
satisfies $$A(w_{1}w_{2},t)=A(w_{1},w_{2}t)A(w_{2},t)$$
which we compute explicitly. These results will not be used
elsewhere in the paper.\\

Let $\sigma:\Bg\rightarrow N(T)$ be a homomorphism making
\eqref{eq:extension} commute. The restriction of $\sigma$
to the pure braid group $\Pg$ maps into $T$ and therefore
factors through the abelianisation of $\Pg$. The following
gives an explicit description of the latter as a $W=\Bg/
\Pg$--module.

\begin{proposition}\label{pr:abelianisation}
Let $Z$ be the free abelian group with one generator
$\gamma_{\alpha}$ for each positive root $\alpha$ of
$\g$ and let $W$ act on $Z$ by $w\medspace \gamma
_{\alpha}=\gamma_{|w\alpha|}$, where $|w\alpha|$ is
equal to $\pm w\alpha$ according to whether $w\alpha$
is positive or negative. Then, 
\begin{enumerate}
\item[(i)] the assignment $\gamma_{\alpha_{i}}\rightarrow
S_{i}^{2}$ extends uniquely to a $W$--equivariant
isomorphism $Z\cong \Pg/[\Pg,\Pg]$.
\item[(ii)] Under the Hurewicz isomorphism $\Pg/[\Pg,\Pg]
\cong H_{1}(\hreg,\IZ)$, $\gamma_{\alpha}$ is mapped onto
a positively oriented simple loop around the hyperplane
$\Ker(\alpha)$.
\end{enumerate}
\end{proposition}
\proof
(i) is proved in \cite[Thm. 2.5]{Ti}. (ii) under the Hurewicz
isomorphism, the action of $W$ on $\Pg/[\Pg,\Pg]$ coincides
with its natural geometric action on $H_{1}(\hreg,\IZ)$. By
Brieskorn's description of the generators of $\Bg$, $\gamma
_{\alpha_{i}}=S_{i}^{2}$ is mapped onto a positively oriented
simple loop around the hyperplane $\Ker(\alpha_{i})$ so that
(ii) follows by $W$--equivariance \halmos\\

For any positive root $\alpha$, pick an element $\lambda
_{\alpha}\in\h$ such that $\exp(2\pi i\lambda_{\alpha})=
\sigma(\gamma_{\alpha})$ and consider the flat connection
on $\hreg\times V$ given by
$$\nabla_{\sigma}=d-\sum_{\alpha\in R^+}
\frac{d\alpha}{\alpha}\cdot\lambda_{\alpha}$$
Fix a basepoint $t_{0}\in\hreg$ and identify $\wt\hreg$
with the space of paths in $\hreg$ pinned at $t_{0}$, modulo
homotopy equivalence. Denote by $\P_{\sigma}(p)\in T$ parallel
transport with respect to $\nabla_{\sigma}$ along one such
path $p$. Then,

\begin{proposition}\label{pr:cocycle}\hfill
\begin{enumerate}
\item[(i)] The map $\wt{\hreg}\times V\rightarrow\hreg
\times V$ given by 
$$(p,v)\rightarrow (p(1),\P_{\sigma}(p)v)$$
descends to an isomorphism of vector bundles $\iota:\wt{\hreg}
\times_{\Pg}V\cong\hreg\times V$.
\item[(ii)] The right action of $\Bg$ on $\wt{\hreg}
\times V$ descends, via $\iota$, to one of $W$ on
$\hreg\times V$ given by
$$w\medspace (t,v)=
(w^{-1}t,\P_{\sigma}(w^{-1}p_{t})\P_{\sigma}(\wt\gamma)
\sigma(\gamma)^{-1}\P_{\sigma}(p_{t})^{-1}v)$$
where $p_{t}$ is any pinned path in $\hreg$ with $p_{t}
(1)=t$, $\gamma\in\Bg$ is any element with image $w$ and
$\wt\gamma$ is its lift to a path in $\hreg$ with $\wt
\gamma(0)=t_{0}$.
\end{enumerate}
\end{proposition}
\proof
One readily checks, by using proposition \ref{pr:abelianisation},
that the monodromy $\Pg\rightarrow T$ of $\nabla_{\sigma}$
coincides with the restriction of $\sigma$ to $\Pg$ from
which (i) follows at once. (ii) is a simple computation
\halmos

\subsection{Tits extensions}\label{ss:Tits}

Let $\sigma:\Bg\rightarrow N(T)$ be a homomorphism making
the diagram \eqref{eq:extension} commute. Tits has given
a simple construction of a canonical, but not exhaustive,
class of such $\sigma$ which differ from each other via
conjugation by an element of $T$. We summarise below the
properties of this class obtained in \cite{Ti}. For any
simple root $\alpha_{i}$, $i=1\ldots n$, let $SL_{2}(\IC)
\cong G_{i}\subseteq G$ be the subgroup with Lie algebra
spanned by $e_{\alpha_{i}},f_{\alpha_{i}},h_{\alpha_{i}}$,
$T_{i}=\exp(\IC\cdot h_{\alpha_{i}})\subset G_{i}$ its
torus and $N_{i}$ the normaliser of $T_{i}$ in $G_{i}$.
Denote by $s_{i}\in W$ the orthogonal reflection corresponding
to $\alpha_{i}$.

\begin{proposition}\label{pr:Tits}\hfill
\begin{enumerate}
\item[(i)] For any choice of $\sigma_{i}\in N_{i}\setminus
T_{i}$, $i=1\ldots n$, the assignment $S_{i}\rightarrow
\sigma_{i}$ extends uniquely to a homomorphism $\sigma:
\Bg\rightarrow N(T)$ making \eqref{eq:extension} commute.
\item[(ii)] If $\sigma,\sigma':\Bg\rightarrow N(T)$ are
the homomorphisms corresponding to the choices $\{\sigma
_{i}\}_{i=1}^{n}$ and $\{\sigma_{i}'\}_{i=1}^{n}$ respectively,
there exists $t\in T$ such that, for any $S\in\Bg$
$$\sigma(S)=t\sigma'(S)t^{-1}$$
\item[(iii)] For any such $\sigma:\Bg\rightarrow N(T)$,
the subgroup $\sigma(\Bg)\subset N(T)$ is an extension
of $W$ by $\IZ_{2}^{n}$ canonically isomorphic to the
group generated by the symbols $a_{i}$, $i=1\ldots n$
subject to the relations
\begin{align}
\underbrace{a_{i}a_{j}\cdots}_{m_{ij}}&=
\underbrace{a_{j}a_{i}\cdots}_{m_{ij}}\label{eq:braid a}\\
a_{i}^{2}a_{j}^{2}&=a_{j}^{2}a_{i}^{2}
\label{eq:braid b}\\
a_{i}^{4}&=1
\label{eq:braid c}\\
a_{i}a_{j}^{2}a_{i}^{-1}&=a_{j}^{2}a_{i}^{-2\<\cor{i},\alpha_{j}\>}
\label{eq:braid d}
\end{align}
for any $1\leq i\neq j\leq n$, where the number
$m_{ij}$ of factors on each side of \eqref{eq:braid a}
is equal to the order of $s_{i}s_{j}$ in $W$.
The isomorphism is given by sending $a_{i}$
to $\sigma_{i}$.
\end{enumerate}
\end{proposition}
\proof
(i) We must show that the $\sigma_{i}$ satisfy
the braid relations \eqref{eq:braid relations}.
For any $1\leq i\neq j\leq n$, set $s_{ij}=s_{i}s_{j}
\cdots\in W$ and $\sigma_{ij}=\sigma_{i}\sigma_{j}
\cdots\in N(T)$ where each product has $m_{ij}-1$ 
factors. The braid relations in $W$ may be written 
as $s_{ij}s_{j'}=s_{j}s_{ij}$ where $j'=j$ or $i$ 
according to whether $m_{ij}$ is even or odd. Thus,
$s_{ij}^{-1}s_{j}s_{ij}=s_{j'}$ and therefore,
$$\delta_{ij}=
\sigma_{j'}^{-1}\sigma_{ij}^{-1}\sigma_{j}\sigma_{ij}
\in T\cap\left(\sigma_{j'}^{-1}\sigma_{ij}^{-1}N_{j}\sigma_{ij}\right)
 = T\cap\sigma_{j'}^{-1}N_{j'}=T_{j'}$$
Repeating the argument with $i$ and $j$ permuted,
we find that $\delta_{ji}\in T_{i'}$ with $i'=i$
or $j$ according to whether $m_{ij}$ is even or
odd. Thus, $\delta_{ij}=\delta_{ji}^{-1}\in T_{i'}
\cap T_{j'}=\{1\}$ where the latter assertion
follows from the simple connectedness of $G$,
and the $\sigma_{i}$ satisfy \eqref{eq:braid relations}.

(ii) Let $t_{i}\in T_{i}$ be such that $\sigma_{i}=
\sigma_{i}'t_{i}$ and choose $c_{i}\in\IC$ such that
$t_{i}=\exp(c_{i}h_{\alpha_{i}})$. Since
$$(s_{i}-1)\medspace\sum_{j=1}^{n}c_{j}\lambda_{j}^{\vee}=
-c_{i}h_{\alpha_{i}}$$
where the $\lambda_{i}^\vee\in\h$ are the fundamental
coweights defined by $\alpha_{i}(\lambda_{j}^{\vee})=
\delta_{ij}$, we find
$$\exp(-\sum_{j}c_{j}\lambda_{j}^{\vee})\sigma_{i}'
\exp(\sum_{j}c_{j}\lambda_{j}^{\vee})
=
\sigma_{i}'\exp(c_{i}h_{\alpha_{i}})
=
\sigma_{i}$$
so that $\sigma$ and $\sigma'$ are conjugate.

(iii) The $\sigma_{i}$ satisfy
\eqref{eq:braid a}--\eqref{eq:braid d} since
$x_{j}^{2}=\exp(i\pi\cor{j})$ for any $x_{j}
\in N_{j}\setminus T_{j}$. Let $K_{\sigma}
\cong\IZ_{2}^{n}$ be the group generated by
the $\sigma_{i}^{2}$ and $K_{\sigma}\subset
\overline{K}_{\sigma}\subset\sigma(\Bg)$ the
kernel of the projection $\sigma(\Bg)\rightarrow
W$. By \eqref{eq:braid d},
$K_{\sigma}$ is a normal subgroup of $\sigma(\Bg)$
and $\sigma(\Bg)/K_{\sigma}$ is generated by
the images $\overline{\sigma_{i}}$ of $\sigma_{i}$
which, in addition to the braid relations
satisfy $\overline{\sigma}_{i}^{2}=1$. Thus,
$\sigma(\Bg)/K_{\sigma}$ is a quotient of $W$,
$K_{\sigma}=\overline{K}_{\sigma}$ and $\sigma
(\Bg)/K_{\sigma}\cong W$. The same argument
shows that if $\Gamma$ is the abstract group
generated by $a_{1},\ldots,a_{n}$ subject to
\eqref{eq:braid a}--\eqref{eq:braid d}, and
$A\subset\Gamma$ is the subgroup generated by
the $a_{i}^{2}$, then $\Gamma/A\cong W\cong\sigma
(\Bg)/K_{\sigma}$. But $A$ is a quotient of
$\IZ^{2}_{n}$ so that the canonical surjection
of $\Gamma$ onto $\sigma(\Bg)$ is an isomorphism
of $A$ onto $K_{\sigma}$ and therefore an isomorphism
of $\Gamma$ onto $\sigma(\Bg)$ \halmos\\

We shall henceforth only use homomorphisms $\sigma:\Bg
\rightarrow N(T)$ of the form given by proposition \ref
{pr:Tits} and refer to them, or their image $\wt{W}=
\sigma(\Bg)\subset N(T)$ as {\it Tits extensions} of
$W$. Note that, given a choice of simple root vectors
$e_{\alpha_{i}},f_{\alpha_{i}}$, $i=1\ldots n$, any
element of $N_{i}\setminus T_{i}$ is necessarily of
the form
\begin{equation}\label{eq:triple exp}
\begin{split}
\sigma_{i}(t_{i})
&=
\exp(t_{i}e_{\alpha_{i}})
\exp(-t_{i}^{-1}f_{\alpha_{i}})
\exp(t_{i}e_{\alpha_{i}})\\
&=
\exp(-t_{i}^{-1}f_{\alpha_{i}})
\exp(t_{i}e_{\alpha_{i}})
\exp(-t_{i}^{-1}f_{\alpha_{i}})
\end{split}
\end{equation}
for a unique $t_{i}\in\IC^{*}$ so that a Tits extension
may be given by choosing elements $t_{1},\ldots,t_{n}
\in\IC^{*}$.

\section
{Generic irreducibility of monodromy representations}
\label{se:gen irred}

\subsection{} In this section, we study in detail the
reducibility of the monodromy of a flat connection of
the form \eqref{eq:kohno}, namely
$$\nabla=d-h\sum_{i\in I}\frac{d\phi_{i}}{\phi_{i}}\cdot r_{i}$$
where the residues $r_{i}$ act on the \fd vector space $V$
and are assumed to satisfy the relations \eqref{eq:relations}.
Let $$\rho_{h}:\pi_{1}(B\setminus\A)\longrightarrow GL(V)$$
be the corresponding one--parameter family of monodromy
representations. If $V$ is reducible under the $r_{i}$,
$\rho_{h}$ is clearly reducible for all values of $h$. We wish
to prove a converse statement. To formulate it, we need the
following\\
 
\definition An analytic curve $\rho_h$ of representations
of a \fg group $\Gamma$ is {\em generically  irreducible}
if $\rho_{h}$ is irreducible for all $h$ in the complement
of an analytic set.\\

The main result of this section is the following
\begin{theorem}\label{th:analytic irred}
If $V$ is irreducible under the $r_{i}$, the monodromy
representation $\rho_{h}$ is generically irreducible.
\end{theorem}

The proof of theorem \ref{th:analytic irred} occupies the
rest of this section; in \S \ref{ss:gen red nablak}, we apply
this result to the monodromy of the connection $\nablak$.
It proceeds by noting that, because reducibility is a closed
condition, an analytic curve $\rho_{h}$ of representations
is either generically irreducible or reducible for all $h$. In
the latter case, we prove the existence of an analytic, {\it
multivalued} curve germ $U(h^{1/m})$ of proper subspaces
of $V$ invariant under $\rho_{h}$. A simple calculation then
shows that the subspace $\left.U(h^{1/m})\right |_{h^{1/m}
=0}\subsetneq V$ is invariant under the $r_{i}$.
Note that in general, a single--valued germ of invariant
subspaces may fail to exist. For example, the curve $c(h)$
of reducible representations of $\Gamma=\IZ$ given by
$$c(h)=\begin{pmatrix} 1& 1 \\ h & 1 \end{pmatrix}$$
with $h\in\IC\setminus\{\pm 1\}$, only admits the multivalued
family of eigenlines $(1,\pm\sqrt{h})$. We will see that such
branching is the worst behaviour that can occur.

\subsection{The variety of reducible representations}
\label{ss:Hom red}
 
Let $\Gamma$ be a \fg group. The set of representations
$\Hom(\Gamma,GL(V))$ has the structure of an affine variety defined
over $\IQ$. Indeed, if $\{\gamma_{1},\ldots,\gamma_{r}\}$ is a system
of generators of $\Gamma$, $\Hom(\Gamma,GL(V))\subset GL(V)^{r}$
is the subset of $r$--tuples of elements satisfying the relations which
define $\Gamma$.
 
\begin{proposition}\label{pr:reducible}
The set $\Hom^{\red}(\Gamma,GL(V))$ of reducible representations
is a Zariski closed subset of $\Hom(\Gamma,GL(V))$.
\end{proposition}
\proof Let $\Gr_{p}(V)$ be the Grassmannian of $p$--planes in $V$
and set
\begin{gather}\label{eq:Rp}
\R_{p}(\Gamma)=
\{(\rho,U)\in\Hom(\Gamma,GL(V))\times\Gr_{p}(V)|\rho(\Gamma)U=U\}
\end{gather}

We claim that $\R_{p}(\Gamma)$ is a Zariski closed subset of
$\Hom(\Gamma,GL(V))\times\Gr_{p}(V)$. Indeed, regarding $U\in
\Gr_{p}(V)$ as a decomposable $p$--tensor $\Lambda=u_{1}\wedge
\cdots\wedge u_{p}\in\IP(\bigwedge^{p}V)$ via the Pl\"{u}cker embedding,
we see that the invariance of $U$ under $\rho\in\Hom(\Gamma,GL(V))$
is equivalent to the relations
\begin{equation}\label{eq:incidence}
\Lambda\wedge\rho(\gamma_{i})\Lambda=0
\end{equation}
for all $i=1\ldots r$, where $\wedge$ is the exterior multiplication in
$\bigwedge^{*}(\bigwedge^{p}V)$. Since the projection $p_{1}:\Hom
(\Gamma,GL(V))\times\Gr_{p}(V)\rightarrow\Hom(\Gamma,GL(V))$
is closed \cite[pg. 76]{Mu}, the set
\begin{equation*}
\begin{split}
\Hom^{\red,p}(\Gamma,GL(V))
&=
\{\rho\in\Hom(\Gamma,GL(V))|
\thinspace\exists\medspace U\in\Gr_{p}(V),
\thinspace\rho\medspace U=U\}\\
&=
p_{1}(\R_{p}(\Gamma))
\end{split}
\end{equation*}
is a closed subset of $\Hom(\Gamma,GL(V))$ and therefore
so is
$$\Hom^{\red}(\Gamma,GL(V))=
\bigcup_{p=1}^{\dim V-1}\Hom^{\red,p}(\Gamma,GL(V))$$
\halmos

\begin{corollary}\label{co:dichotomy}
An analytic curve $\rho_{h}:\Gamma\longrightarrow GL(V)$ of
representations is either generically irreducible or reducible
for all values of $h$.
\end{corollary}

\remark The variety $\R^p(\Gamma)$ is clearly defined over $\IQ$.
By \cite[prop. 1.1.1]{St}, this is also true of $\Hom^{\red,p}(\Gamma,
GL(V))$ since it is defined over $\ol{\IQ}$ by \cite[pg. 76]{Mu} and
the underlying set of $\ol{\IQ}$--points is stable under the action of
$\Gal(\ol{\IQ},\IQ)$. Finally, the regular map
$$\pi_p:\R_p(\Gamma)
\longrightarrow\Hom^{\red,p}(\Gamma,GL(V))$$
induced by projection on the first factor is also defined over $\IQ$
since $p_1$ is.

\subsection{The existence of a multivalued section}

Retain the notation of \S \ref{ss:Hom red}. We will be concerned in this
subsection with constructing a multivalued section to the projection $\pi
_{p}$ over a curve germ in $\Hom^{\red,p}(\Gamma,GL(V))$. Let
$\Gamma_r$ be the free group $\Gamma_r$ on $r$ generators.
The commutative diagram 
$$\begin{CD}
\R^{p}(\Gamma)	      @>>> \R^{p}(\Gamma_r)	\\
    @V{\pi_p}VV           	@V{\pi_p}VV	\\ 
\Hom^{\red,p}(\Gamma) @>>> \Hom^{\red,p}(\Gamma_r)
\end{CD}$$
shows that it suffices to find such a section for $\Gamma_r$.
For the remainder of this subsection we therefore assume
that $\Gamma=\Gamma_{r}$.\\

We shall in fact prove the existence of a {\it formal} multivalued
section. The same argument yields an analytic multivalued one.
We need some notation. Let $A$ be the coordinate ring of $\Hom
(\Gamma,GL(V))\cong GL(V)^{r}$, $\a_p,\a\subset A$ the ideals
of $\Hom^{\red,p}(\Gamma,GL(V))$, $\Hom^{\red}(\Gamma,GL(V))$
respectively and $S_p=A/\a_p$, $S=A/\a$ their coordinate rings.
Let $\wt{\Gr}_p(V)$ be the affine cone in $\bigwedge^p(V)$ defined
by the Pl\"ucker equations, $B_{p}=\IC[\Hom(\Gamma,GL(V))\times
\wt{\Gr}_{p}(V)]$ and $\b_p\subset B_{p}$ the ideal defined by
equations \eqref{eq:incidence}. Set
$$R_{p}=B_p/\b_p=\IC[\wt{\R}^{p}(\Gamma)]$$
where $\wt{\R}^{p}(\Gamma)\subset\Hom(\Gamma,GL(V))\times
\wt{\Gr}_{p}(V)$ is the preimage of the variety $\R_{p}(\Gamma)$
defined by \eqref{eq:Rp}. The projection $\pi_p$ induces a ring 
homomorphism $\psi_p:S_{p}\to R_{p}$.\\

Let $\rho$ be a fixed reducible representation. A formal curve
of representations starting at $\rho$ is a homomorphism of $\IC
$--algebras
$$\phi:A\longrightarrow\IC[[h]]$$
where $\IC[[h]]$ is the ring of formal power series in $h$
whose reduction mod $h$ is the $\IC$--point $\rho$. $\phi$ is a
curve of reducible representations if it descends to $S$. Let $\p$
be the kernel of $\phi$. $\p$ is a prime ideal since $\IC[[h]]$ is an
integral domain and $\p\supset\a=\bigcap_{p=1}^{n-1}\a_p$ so
that $\p\supset\a_p$ for some $p$ and we obtain a homomorphism
$\phi:S_p\longrightarrow\IC[[h]]$.\\

\definition A multivalued formal section of $\pi_p$ is a homomorphism
$\chi$ fitting into the following
commutative diagram
$$\begin{CD} 
\IC[[k]] \ @<\chi<<R_p\\
@A{f_m}AA     @A{\psi_p}AA \\ 
\IC[[h]]@<{\phi}<<S_p
\end{CD}$$
where $f_m(h) = k^m$.\\

Let $\F=\IC((h))$ be the field of fractions of $\IC[[h]]$, $\ol{\F}$ its
algebraic closure and $\iota:\F\to\ol{\F}$ the corresponding inclusion.
The key step in finding a formal multivalued section is the following

\begin{proposition}\label{pr:completeness}
There exists a homomorphism of $\IC$--algebras $\tau:R_p
\to\ol{\F}$ such that the following diagram is commutative
\begin{equation}\label{eq:multi}
 \begin{CD} 
\ol{\F} \ @<\tau<<R_p\\
@A{\iota}AA     @A{\psi_p}AA \\ 
\IC[[h]]@<{\phi}<<S_p
\end{CD}
\end{equation}
and $\tau(x_{I})\neq 0$ for some $I$, where $x_I$ are the Pl\"ucker
coordinates on $\wt{\Gr}_p(V)$.
\end{proposition}
\proof Note first that $\ol{\F}$ is isomorphic to $\IC$ since both are
algebraically closed fields of uncountable transcendence degree
over $\IQ$, and choose an isomorphism $\phi:\IC\rightarrow\ol{\F}$.
Since $X=\Hom^{\red,p}(\Gamma,GL(V))$, $Y=\wt{\R}_{p}$ are
defined over $\IQ$, $\phi$ induces an isomorphism between the
the $\IC$ and $\ol{\F}$--points of $X$ (resp. $Y$) which we denote
by $\ul{X}(\IC),\ul{X}(\ol{\F})$ (resp. $\ul{Y}(\IC),\ul{Y}(\ol{\F})$). Since
the projection $\pi_p$ is defined over $\IQ$, we therefore get a
commutative diagram
$$\begin{CD}
\ul{Y}(\IC)&@>>>&\ul{Y}(\ol{\F})\\
@VVV&      &@VVV\\
\ul{X}(\IC)&@>>>&\ul{X}(\ol{\F})
\end{CD}$$
and the rightmost vertical arrow is onto \halmos\\

We now construct the desired multivalued formal section. Recall that
the algebraic closure of $\IC((h))$ is the field $\Pus{h}$ obtained by
adjoining all the roots of $h$ to $\IC((h))$ \cite[thm. 3.1]{Wal}.
 
\begin{theorem}
There exists an $m\geq 1$ and a homomorphism of $\IC$--algebras
$\chi:\wh{R}_p\to\IC[[k]]$ such that the following diagram is
commutative
$$\begin{CD} 
\IC[[k]] \ @<\chi<<R_p\\
@Af_mAA     @A{\psi_p}AA \\ 
\IC[[h]]@<{\phi}<<S_p
\end{CD}$$
and $\chi(x_{I})\neq 0$ for some $I$.
\end{theorem}
\proof Proposition \ref{pr:completeness} yields a homomorphism
$\tau$ taking values in the quotient field $\IC((k))$ of $\IC[[k]]$
for some root $k$ of $h$. Such a homomorphism amounts to
assigning an element of $\IC((k))$ to each of the generators
$g_{ij}^k$, $1\leq i,j\leq\dim(V)$, $1\leq k\leq r$ and $x_I$ of
$R_p$ in such a way that these elements satisfy the
defining equations of $\wt{\R}_p(\Gamma)$. By \eqref{eq:multi},
$\tau(g_{ij}^k)=\iota\phi(g_{ij}^k)\in\IC[[k]]$. Since the equations
are homogeneous in the $x_I$, we can multiple the $\tau(x_I)$
by an appropriate power of $k$ so that the resulting elements
of $\IC((k))$ are in $\IC[[k]]$ \halmos\\

The following result explains the meaning of the $\IC[[k]]
$--point $\chi$.

\begin{proposition}\label{pr:formalpoints}
The homomorphism $\chi$ gives rise to a rank $p$ summand
$\U$ of the module $V[[k]]=V\otimes\IC[[k]]$ together with an
$r$--tuple of invertible (over $\IC[[k]]$) matrices  $g^1(k^m),
\ldots,g^r(k^m)$ leaving $\U$ invariant.
\end{proposition}

\subsection{Proof of theorem \ref{th:analytic irred}}

Let $\Gamma$ be the fundamental group $\pi_{1}(B\setminus
\A)$. It is well--known that $\Gamma$ is finitely--generated,
see \eg \cite[prop. A2, pg. 181]{BMR}. If the curve $\rho
_{h}$ of monodromy representations is not generically
irreducible, it lies, by corollary \ref{co:dichotomy}, in the
variety of reducible representations and therefore in some
$\Hom^{\red,p}(\Gamma,GL(V))$, $1\leq p\leq\dim
V-1$. Denote by $\IC\{k\}$ the ring of convergent power
series at $k=0$. Let $\chi$ be a multivalued analytic section
and let $\U$ be the corresponding rank $p$ summand of the
free $\IC\{k\}$--module $\V=V\otimes\IC\{k\}$ obtained by
applying the analytic version of proposition \ref
{pr:formalpoints}. $\U$ corresponds to an multivalued 
analytic curve $U(k)$ in $\Gr_p(V)$ which is invariant under
$\rho_h$.\\

Set $\wt{\rho}_k=\rho_{k^m}$ and let $\wt{\rho}$ be the
element of $\Hom(\Gamma,GL(\V))$ corresponding to $\wt
{\rho}_k$. Since $\wt{\rho}_k$ is an analytic function
of $h=k^m$, the elements 
$$A(\gamma)=\frac{I-\wt{\rho}(\gamma)}{k^m}$$
leave $\V$ invariant. It follows that $A(\gamma)\medspace
\U\subset\U$ since $\wt{\rho}(\Gamma)\medspace\U=\U$. The
following standard result shows that the subspace $U(0)
\subsetneq V$ is invariant under the residues $r_{i}$,
thus concluding the proof of theorem \ref{th:analytic irred} \halmos

\begin{lemma}\label{le:generator}
For each hyperplane $H_{j}$ of the arrangement $\A$, there
exists an element $\gamma_{j}\in\pi_{1}(B\setminus\A)$ such
that
$$\rho_{h}(\gamma_{j})=1+2\pi ih\cdot r_{j} \mod h^{2}$$
\end{lemma}
\proof Let $x_{0}\in B\setminus\A$ be a base point and 
$\gamma_{j}\in\pi_{1}(B\setminus\A;x_{0})$ a generator
of monodromy around $H_{j}$ (see, \eg \cite[pg. 180--1]{BMR}).
Recall that such an element is obtained as follows. Choose a
path $p:[0,1]\to B$ such that
$$p(0)=x_{0},\quad
p([0,1[)\subset B\setminus \A
\quad\text{and}\quad
p(1)\in H_{j}\setminus\bigcup_{j'\neq j}H_{j'}$$
Let $D$ be a small ball centred at $p(1)$ and contained in
$B\setminus\bigcup_{j'\neq j}H_{j'}$, let $u\in[0,1[$ be
such that $p(s)\in D$ for any $s\geq u$ and let $\ell$ be
a positively oriented generator of $\pi_{1}(D\setminus
H_{j};p(u))\cong\IZ$. Then,
$$\gamma_{i}=p_{u}^{-1}\cdot\ell\cdot p_{u}$$
where $p_{u}(t)=p(ut)$ and the concatenation of paths
is read from right to left. Picard iteration readily
yields that, mod $h^{2}$,
$$\rho_{h}(\gamma_{j})
=
1+h\sum_{i\in I}
\int_{\gamma_{j}}\frac{d\phi_{i}}{\phi_{i}}\cdot r_{i}
=
1+h\sum_{i\in I}
\int_{\ell}\frac{d\phi_{i}}{\phi_{i}}\cdot r_{i}
=
1+2\pi ih\cdot r_{j}$$
where the last equality follows from the residue theorem
since the forms $d\phi_{j'}/\phi_{j'}$, $j'\neq j$ do not
have any poles in $D$ \halmos\\

The above proof readily yields the following analogue
of theorem \ref{th:analytic irred} in the formal setting. Let
$$\wt{\rho}:\pi_{1}(B\setminus\A)\longrightarrow GL(V\La)$$
be the representation obtained by regarding $\rho_{h}$
as formal in $h$, letting $\pi_{1}(B\setminus\A)$ act
on $V[[h]]$ and extending coefficients to $V((h))$. Then

\begin{theorem}\label{th:formal irred}
If $V$ is irreducible under the $r_{i}$, $\wt{\rho}$ is
absolutely irreducible, that is irreducible over the
algebraic closure $\Pus{h}$ of $\IC\La$.
\end{theorem}

\subsection{Generic irreducibility of the monodromy of $\nablak$}
\label{ss:gen red nablak}

Assume now that $V$ is a $\g$--module and let
$\wt{W}=\sigma(\Bg)$ be a Tits extension with
sign group $\Sigma=\sigma(\Pg)\cong\IZ_{2}^{n}$.
Let 
$$\rho_{h}^{\sigma}:\Bg\longrightarrow GL(V)$$
be the corresponding one--parameter family of
monodromy representations defined by proposition
\ref{th:existence} and $\rho^{\sigma}:\Bg\longrightarrow
GL(V\La)$ the representation obtained by regarding
$\rho_{h}^{\sigma}$ as formal in $h$ and extending
coefficients to $\IC\La$.

\begin{theorem}\label{th:Bg irred}
Let $U\subseteq V$ be a subspace invariant, and
irreducible under the Casimirs $\kalpha$ and
$\wt{W}$. Then,
\begin{enumerate}
\item $\rho_{h}^{\sigma}:\Bg\longrightarrow GL(U)$
is generically irreducible.
\item $\rho^{\sigma}:\Bg\longrightarrow GL(U\La)$
is absolutely irreducible.
\end{enumerate}
\end{theorem}
\proof
(i) Assume $\rho_{h}^{\sigma}$ to be reducible for
all $h$. Proceeding as in the proof of theorem \ref
{th:analytic irred}, we find a $1\leq p\leq\dim U-1$
and a rank $p$--summand $\X$ of $\U=U\otimes\IC\{k\}$
invariant under the germ of $\rho_{h}^{\sigma}$ at
$h=0$. In particular, $\X(0)$ is invariant under
$\rho_{0}^{\sigma}(\Bg)=\wt{W}$. For any positive
root $\alpha$, let $\gamma_{\alpha}\in\Pg$ be the
generator of monodromy around the hyperplane $\Ker
(\alpha)$ given by lemma \ref{le:generator}. Note
that $\sigma(\gamma_{\alpha})$ lies in the sign
group $\Sigma$ and therefore has order $1$ or $2$.
By lemma \ref{le:generator} and proposition \ref
{pr:recipe} we find that, mod $h^{2}$,
$$\rho_{h}^{\sigma}(\gamma_{\alpha}^{2})=
\sigma(\gamma_{\alpha})^{2}(1+2\pi ih\cdot\kalpha)^{2}=
1+4\pi ih\cdot\kalpha$$
so that $\X(0)$ is also invariant under the
Casimirs $\kalpha$. The proof of (ii) is identical
\halmos\\

Similarly, we obtain

\begin{theorem}\label{th:Pg irred}
Let $U\subseteq V$ be a subspace invariant, and
irreducible under the Casimirs $\kalpha$ and
$\Sigma$. Then,
\begin{enumerate}
\item $\rho_{h}^{\sigma}:\Pg\longrightarrow GL(U)$
is generically irreducible.
\item $\rho^{\sigma}:\Pg\longrightarrow GL(U\La)$ 
is absolutely irreducible.
\end{enumerate}
\end{theorem}

We specialise our results further to the case where $U$ is
the zero weight space $V[0]$ of $V$. Recall that the latter
is canonically acted upon by $W\cong N(T)/T$ so that the
restriction of $\rho_{h}^{\sigma}$ to $V[0]$ does not depend
upon the choice of $\sigma$. We owe the following somewhat
surprising observation to B. Kostant

\begin{proposition}\label{pr:basic criterion}
$V[0]$ is irreducible under the Casimirs $\kalpha$
iff it is irreducible under the $\kalpha$ and $W$.
\end{proposition}
\proof
The simple reflection $s_{i}\in W$ acts on the zero
weight space $V^{i}_{n}[0]$ of the irreducible
$\sl{2}^{\alpha_{i}}$--module of dimension $2n+1$
as multiplication by $(-1)^{n}$. Thus, if $p^{i}_
{\varepsilon}$, $\varepsilon=0,1$ are the spectral
projections for the restriction of $C_{\alpha_{i}}$
to $V[0]$ corresponding to the Casimir eigenvalues
of $V^{i}_{n}$, with $n=\epsilon\mod 2$, $s_{i}$
acts on $V[0]$ as $p^{i}_{0}-p^{i}_{1}$ and is
therefore a polynomial in $C_{\alpha_{i}}$. It
follows that a subspace $U\subseteq V[0]$ invariant
under the $\kalpha$ is also invariant under $W$
\halmos

\begin{corollary}\label{co:gen red V[0]}
The following statements are equivalent
\begin{enumerate}
\item $V[0]$ is irreducible under the Casimirs $\kalpha$.
\item $V[0]$ is generically irreducible under $\Pg$.
\item $V[0]\La$ is absolutely irreducible under $\Pg$.
\item $V[0]$ is generically irreducible under $\Bg$.
\item $V[0]\La$ is absolutely irreducible under $\Bg$.
\end{enumerate}
\end{corollary}

\section{The Casimir algebra $\Cg$ of $\g$}\label{se:Casimir}

Recall from the Introduction that the {\it Casimir algebra}
$\Cg$ of $\g$ is the algebra
$$\Cg=\<\kalpha\>_{\alpha\in R_{+}}\vee\h\subseteq\Ugh$$
generated by the $\kalpha$, or equivalently the Casimirs
$\calpha$, and $\h$ inside $\Ugh$.

\subsection{The Casimir algebra of $\sl{3}$}\label{ss:sl3}

\begin{theorem}\label{th:Csl3}
If $\g=\sl{3}$, then
\begin{enumerate}
\item $\Cg$ is a proper subalgebra of $\Ugh$.
\item $\Ugh$ is generated by  $\Cg$ and the center of $\Ug$.
In particular, the Casimirs $\kalpha$ act irreducibly on the
weight spaces of simple $\g$--modules.
\end{enumerate}
\end{theorem}

We need some preliminary results. Let $e_{1},\ldots,e_{n}$ be the
canonical basis of $\IC^{n}$ and $E_{ij}e_{k}=\delta_{jk}e_{i}$ the
corresponding elementary matrices. For any sequence $I=(i_1,
\ldots,i_k)$ of distinct elements of $\{1,\ldots,n\}$, set
$$E_I=E_{i_1i_2}E_{i_2i_3}\cdots E_{i_{k-1}i_k}E_{i_ki_1}
\in U\sl{n}^{\h}$$

\begin{proposition}\label{pr:Sgh}
If $\g=\sl{n}$, $\Ugh$ is generated as an algebra by $\h$ and the
monomials $E_I$ corresponding to sequences $(i_{1},\ldots,i_{k})$
such that $i_{1}=\min_{l}i_{l}$.
\end{proposition}
\proof It suffices to show that the $E_{I}$ and $\h$ generate $\gr(\Ug^\h)
=S\g^\h$. $\Sgh$ is clearly spanned by elements of the form $p\cdot
E_{I,\sigma}$ where $p\in S\h$, $I=(i_{1},\ldots,i_{k})$ is a sequence
of elements in $\{1,\ldots,n\}$, $\sigma\in\SS{k}$ is a permutation
and $E_{I,\sigma}=E_{i_{1}i_{\sigma(1)}}\cdots E_{i_{k}i_{\sigma(k)}}$.
Writing $\sigma$ as a product $\tau_{1}\cdots\tau_{r}$ of disjoint cycles
with $\tau_{j}=(m_{j}^{1}\cdots m_{j}^{k_{j}})$ and setting $I_{j}=(i_{m_
{j}^{1}},\cdots,i_{m_{j}^{k_{j}}})$ shows that, in $S\g$, $E_{I,\sigma}
=E_{I_{1}}\cdots E_{I_{r}}$ \halmos

\begin{corollary}\label{co:Usl3}
If $\g=\sl{3}$, $\Ugh$ is generated by $\h$,
\begin{xalignat*}{3}
F_{12}&=E_{12}E_{21}&
F_{13}&=E_{13}E_{31}&
F_{23}&=E_{23}E_{32}
\end{xalignat*}
and
\begin{xalignat*}{2}
G_{123}&=E_{12}E_{23}E_{31}&
G_{132}&=E_{13}E_{32}E_{21}
\end{xalignat*}
\end{corollary}

{\sc Proof of Theorem \ref{th:Csl3}}. (ii) Let $\theta_{i}-\theta_{j}$,
$1\leq i<j\leq 3$ be the positive roots of $\g$. The $\sl{2}$--triple
corresponding to $\theta_{i}-\theta_{j}$ is $\{E_{ij},E_{ji},E_{ii}-E
_{jj}\}$ and
$$\kappa_{\theta_{i}-\theta_{j}}=
E_{ij}E_{ji}+E_{ji}E_{ij}=2F_{ij}\mod\h$$
$\Sgh$ is therefore generated, as a $\Gr(\Cg)$--algebra,
by $G=G_{123}+G_{132}$ since
$$G_{123}-G_{132}=F_{13}+[F_{23},F_{12}]\in\Cg$$
It therefore suffices to show that $G$ lies in the
algebra generated by $\Cg$ and $Z(\Ug)$, which is a
simple exercise.
(i) Let $\sigma$ be the involution of $\g$ given
by $\sigma(X)=-X^{t}$. $\sigma$ leaves $\h$ invariant and
descends to an involution of $Q=\Ugh/\Ugh\h$ fixing $\Cg$.
It therefore suffices to show that $\sigma$ does not act
trivially on $Q$ or the associated graded $\gr(Q)$.  
However the image of $G=G_{123}+G_{132}$ in $\gr(Q)$
satisfies $\sigma(G)=-G$ since $\sigma(G_{123})=-G_{132}$
in $\gr(Q)$ \halmos

\begin{corollary}
If $V$ be a simple $\sl{3}$--module, the monodromy
of $\nablak$ yields generically irreducible representations of
Artin's pure braid group $P_{3}$ on the weight spaces of $V$
and of the braid group $B_{3}$ on the zero weight space of
$V$.
\end{corollary}

\subsection{The Casimir algebra of $\g\ncong\sl{n}$}
\label{ss:not sln}

Assume that $\g$ is simple and not isomorphic to $\sl{n}$
and let $V$ be the kernel of the commutator map $[\cdot,
\cdot]:\g\wedge\g\rightarrow\g$. It is known that $V$ is
a simple $\g$--module \cite{Re}.

\begin{theorem}\label{th:Casimir not sln}
The zero weight space $V[0]$ is reducible under $\Cg$. In
particular, $\Cg$ and $Z(\Ug)$ do not generate $\Ugh$.
\end{theorem}
\proof Since $\g\wedge\g\cong\g\oplus V$ and the zero
weight space of $\g\wedge\g$ has a basis given by $h_
{i}\wedge h_{j}$, $1\leq i<j\leq n$, and $e_{\alpha}
\wedge f_{\alpha}$, $\alpha\in R_{+}$, where $h_{1},
\ldots,h_{n}$ is a basis of $\h$, we find that
$$\dim V[0]=
\frac{n(n-1)}{2}+|R_{+}|=
\frac{n(n-1)}{2}+\frac{m-3n}{2}$$
where $m=\dim(\g)>3n$. Thus, $\h\wedge\h$ is a proper
subspace of the zero weight space of $V$ and it suffices
to show that it is invariant under the $\kalpha$. This
follows at once from the fact that, for any $t_{1},t_{2}
\in\h$,
\begin{equation*}
\begin{split}
e_{\alpha}f_{\alpha}\medspace t_{1}\wedge t_{2}
&=
e_{\alpha}\left(
\alpha(t_{1})f_{\alpha}\wedge t_{2}+
\alpha(t_{2})t_{1}\wedge f_{\alpha}
\right)\\
&=
\alpha(t_{1})h_{\alpha}\wedge t_{2}+
\alpha(t_{2})t_{1}\wedge h_{\alpha}
\end{split}
\end{equation*}
\halmos

\subsection{A general reducibility criterion for $V[0]$}
\label{ss:Chevalley}

Let $\Theta$ be the Chevalley involution of $\g$
relative to a choice of simple root vectors $e_{
\alpha_{i}}, f_{\alpha_{i}}$, \ie the automorphism
of $\g$ defined by
$$\Theta(e_{\alpha_{i}})=-f_{\alpha_{i}},
\qquad
\Theta(f_{\alpha_{i}})=-e_{\alpha_{i}}
\qquad\text{and}\qquad
\Theta(h_{\alpha_{i}})=-h_{\alpha_{i}}$$

If $V$ is a simple, finite--dimensional $\g$--module,
and $V^{\Theta}$ is the module obtained by twisting
the action of $\g$ by $\Theta$, then $V^{\Theta}$ is
isomorphic to the dual $V^{*}$ of $V$. In particular,
if $V$ is self--dual, there exists an involution
$\Theta_{V}$ acting on $V$ such that, for any $X\in\g$,
$$\Theta_{V}X\Theta_{V}=\Theta(X)$$

Although $\Theta_{V}$ is only unique up to a sign, we
shall abusively refer to it as {\it the} Chevalley
involution of $V$. Since $\Theta$ acts as $-1$ on
the Cartan subalgebra $\h$ and fixes the Casimirs
$\kalpha$, $\Theta_{V}$ leaves the zero weight space
$V[0]$ invariant and commutes with the action of $\Cg$.
The following gives a useful criterion to show that
$\Theta_{V}$ does not act as a scalar on $V[0]$ and
therefore that the latter is reducible under $\Cg$.  

\begin{proposition}\label{pr:Chevalley}
Let $V$ be a self--dual $\g$--module with $V[0]\neq 0$.
Let $\r\subset\g$ be a reductive subalgebra normalised
by $\h$. Assume that there exists a non--zero vector
$v\in V[0]$ lying in a simple $\r$--module $U$ such
that $U\ncong U^*$. Then, $\Theta_{V}$ does not act
as a scalar on $V[0]$ and the latter is reducible under
$\Cg$.
\end{proposition}
\proof The assumptions imply that $\Theta$ leaves $\r$
invariant and therefore acts as a Chevalley involution
on it. Thus, $\Theta_{V}U\subset V$ is a simple $\r
$--module isomorphic to $U^{*}$ which has zero intersection
with $U$ since $U\ncong U^{*}$. In particular, $\Theta
_{V}v$ is not proportional to $v$ \halmos\\

We record for later use the following alternative proof
of theorem \ref{th:Casimir not sln}.

\begin{proposition}\label{pr:ext2g[0]}
Let $\g\ncong\sl{n}$ and let $V$ be the simple, self--dual
$\g$--module $V=\Ker[\cdot,\cdot]\subset\bigwedge^{2}\g$.
Then, the Chevalley involution $\Theta$ does not act as a
scalar on $V[0]$.
\end{proposition}
\proof $\Theta$ acts as $+1$ on the subspace $\h\wedge\h
\subset V[0]$ and as $-1$ on the span of the vectors $e_
{\alpha}\wedge f_{\alpha}\in\bigwedge^{2}\g[0]$. The
conclusion follows since, as noted in the proof of
theorem \ref{th:Casimir not sln}, $\h\wedge\h$ is a
proper subspace of $V[0]$ \halmos

\subsection{The \GZ branching rules}\label{ss:GZ}

Let $e_{1},\ldots,e_{n}$ be the canonical basis of $\IC^{n}$
and $E_{ab}e_{c}=\delta_{bc}e_{a}$ the corresponding elementary
matrices. Consider the chain of subalgebras
\begin{equation}\label{eq:chain}
\gl{1}\subset\gl{2}\subset\cdots\subset\gl{n-1}\subset\gl{n}
\end{equation}
where each $\gl{k}$ is spanned by the matrices $E_{ij}$,
$1\leq i,j\leq k$. By the Gelfand--Zetlin branching rules
\cite{GZ,Zh1}, the irreducible representation $V_{\lambda}$
of $\gl{k}$ with highest weight $\lambda=(\lambda_{1},\ldots,
\lambda_{k})\in\IZ^{k}$ decomposes under $\gl{k-1}$ as
$$\res_{\gl{k}}^{\gl{k-1}}V_{\lambda}=
\bigoplus_{\ol{\lambda}} V_{\ol{\lambda}}$$
where $V_{\ol{\lambda}}$ is the irreducible $\gl{k-1}$--module
with highest weight $\ol{\lambda}$ and $\ol{\lambda}=(\ol{\lambda}
_{1},\ldots,\ol{\lambda}_{k-1})\in\IZ^{k-1}$ ranges over all
dominant weights of $\gl{k-1}$ satisfying the inequalities
$$\lambda_{1}\geq\ol{\lambda}_{1}\geq\lambda_{2}\geq\cdots\geq
\ol{\lambda}_{k-1}\geq\lambda_{k}$$
which we denote by $\lambda\succ\ol{\lambda}$. Since the above
decomposition is multiplicity--free, it follows, by restricting
in stages from $\gl{n}$ to $\gl{1}$ along \eqref{eq:chain}, that
any simple $\gl{n}$--module $V$ possesses a basis labelled by
\GZ patterns, \ie arrays $\mu$ of the form

\newarrow{Dot}{}{.}{}{.}{} 
$$\begin{diagram}[height=1.3em,width=1em]
\mmu{n}{1}&           &\mmu{n}{2}&          & \cdots   &          &\mmu{n}{n-1}  &              &\mmu{n}{n}\\
         &\mmu{n-1}{1}&          &          & \cdots   &          &              &\mmu{n-1}{n-1}&          \\
         &            &		 &          &          &          &		 &              &          \\
         &            &          &\mmu{2}{1}&          &\mmu{2}{2}&              &              &          \\
         &            &          &          &\mmu{1}{1}&          &              &              &
\end{diagram}$$
where the top row $\mmu{n}{}$ is equal to the highest weight of $V$
and each pair $\mmu{k}{}\in\IZ^{k},\mmu{k-1}{}\IZ^{k-1}$ of consecutive
rows satisfies $\mmu{k}{}\succ\mmu{k-1}{}$. Up to a scalar factor,
the vector $v_{\mu}$ corresponding to the above pattern is uniquely
determined by the requirement that it tranforms under each $\gl{k}
\subset\gl{n}$ according to the irreducible representation with
highest weight $\mmu{k}{}$. In particular, since the central element
$\sum_{i=1}^{k}E_{ii}\in\gl{k}$ acts in the latter as multiplication
by $|\mmu{k}{}|=\sum_{i=1}^{k}\mmu{k}{i}$, we find that, for any
$1\leq i\leq n$,
$$E_{ii}v_{\mu}=
\left(|\mmu{i}{}|-|\mmu{i-1}{}|\right)v_{\mu}$$
so that $v_{\mu}$ has weight zero for the action of $\sl{n}$ iff,
for any $1\leq i\leq n$,
$$|\mmu{i}{}|=i|\mmu{1}{}|=i\mmu{1}{1}$$

For later use in \S \ref{ss:irred A}, we shall need the non--vanishing
of some of the matrix coefficients for the action of the simple roots
vectors of $\gl{n}$ in the above basis. This follows from the explicit
formulae for the action of all elementary matrices $E_{ij}$ in
a suitably normalised \GZ basis $v_{\mu}$ which may be found in
\cite{GZ,Zh2}.

\begin{theorem}[Gelfand--Zetlin]\label{th:nonzero}
Let $\mu$ be a \GZ pattern, then, for any $1\leq i\leq n-1$
\begin{align*}
E_{i\medspace i+1} v_{\mu}&=
\sum_{\mu'}c_{\mu,\mu'}^{i} v_{\mu'}\\
\intertext{where the sum ranges over all patterns $\mu'$ obtained
from $\mu$ by adding 1 to one of the entries of its $i$th row and
the coefficients $c_{\mu,\mu'}^{i}$ are non--zero and,}
E_{i+1\medspace i} v_{\mu}&=
\sum_{\mu'}\wt{c}_{\mu,\mu'}^{i} v_{\mu'}
\end{align*}
where the sum ranges over all patterns $\mu'$ obtained from $\mu$
by substracting 1 to one of the entries of its $i$th row and the
coefficients $\wt{c}_{\mu,\mu'}^{i}$ are non--zero.
\end{theorem}

\begin{corollary}\label{co:nonzero}
Let $\mu$ be a \GZ pattern, then
$$E_{i\medspace i+1}E_{i+1\medspace i}v_{\mu}=
\sum_{\mu'}d^{i}_{\mu,\mu'}v_{\mu'}$$
where the sum ranges over all patterns $\mu'$ differing from $\mu$
by the addition and the substraction of 1 on a pair of (not
necessarily distinct) entries of the $i$th row and $d^{i}_
{\mu,\mu'}\neq 0$ if $\mu\neq\mu'$.
\end{corollary}

\subsection{The Casimir algebra of $\sl{n}$, $n\geq 4$}
\label{ss:Casimir sln}

\begin{theorem}\label{th:Casimir sln}
If $\g=\sl{n}$, $n\geq 4$, there exists a simple, self--dual $\g
$--module $V$ such that the Chevalley involution $\Theta$
does not act as a scalar on $V[0]$. Thus, $V[0]$ is reducible
under $\Cg$ and $\Cg$ and $Z(U\g)$ do not generate $\Ugh$.
\end{theorem}
\proof
By proposition \ref{pr:Chevalley}, it suffices to exhibit
an irreducible representation $V$ of $\gl{n}$ which is
self--dual as $\sl{n}$--module and a \GZ pattern $\mu$
describing a zero--weight vector for $\sl{n}$ in $V$
such that, for some $2\leq k\leq n-1$, the $\sl{k}$--module
$U$ with highest weight $\mu^{(k)}=(\mmu{k}{1},\ldots,\mmu
{k}{k})$ is not self--dual. Since the highest weight of
$U^{*}$ is $(-\mmu{k}{k},\ldots,-\mmu{k}{1})$, such a $U$
is self--dual iff the sum $\mmu{k}{i}+\mmu{k}{k+1-i}$ does
not depend upon $i=1\ldots k$. The following is a suitable
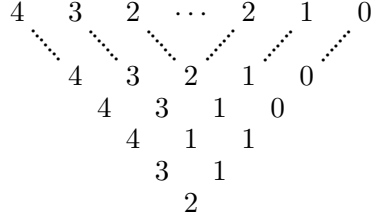
\GZ pattern $\mu$

$$\begin{diagram}[height=1.1em,width=1em]
4&      &3&      &2&      &\cdots&      &2&      &1&      &0 \\
 &\rdDot& &\rdDot& &\rdDot&      &\ldDot& &\ldDot& &\ldDot&  \\
 &      &4&      &3&      &2     &      &1&      &0&      &  \\
 &      & &4     & &3     &      &1     & &0     & &      &  \\
 &      & &      &4&      &1     &      &1&      & &      &  \\
 &      & &      & &3     &      &1     & &      & &      &  \\
 &      & &      & &      &2     &      & &      & &      &
\end{diagram}$$

since the $\sl{3}$--module with highest weight $(4,1,1)$ is not
self--dual and, for any $n\geq 4$, the $\sl{n}$--module with
highest weight
$$(4,3,\underbrace{2,\ldots,2}_{n-4},1,0)$$
is self--dual \halmos

\section{A conjecture of Kwon and Lusztig on quantum Weyl groups}
\label{se:Kwon}

\subsection{} We discuss below some results of Kwon on
$q$--Weyl group actions of Artin's braid group $B_{n}$
on the zero weight spaces of $\Uhsl{n}$--modules \cite
{Kw}. We disprove in particular a conjecture of his and
Lusztig's stating the irreducibility of all such
representations.\\

Let $\Uhg$ be the Drinfeld--Jimbo quantum group corresponding
to $\g$ \cite{Dr1,Ji}, which we regard as a Hopf algebra over the
ring $\ICh$ of formal power series in the variable $\hbar$. By a
{\it finite--dimensional} representation of $\Uhg$ we shall mean
a $\Uhg$--module $\V$ which is topologically free and \fg over
$\ICh$. The isomorphism class of such a representation is uniquely
determined by that of the $\g$--module $V=\V/\hbar\V$ \cite{Dr2}.\\

Lusztig, and independently \KR and Soibelman \cite{Lu,KR,So},
proved that any such $\V$ carries an action, called the {\it
quantum Weyl group action} of the braid group $\Bg$.
Its reduction mod $\hbar$ factors through the Tits extension
$\wt W$ given by the triple exponentials \eqref{eq:triple exp}
with $t_{i}=1$. Specifically, this action is given by mapping
the generator $S_{i}$ of $\Bg$ to the triple $q$--exponential
\cite{Sa}
\begin{equation}\label{eq:triple qexp}
\begin{split}
&
\exp_{q_{i}^{-1}}(q_{i}^{-1}E_{i}q_{i}^{-H_{i}})
\exp_{q_{i}^{-1}}(-F_{i})
\exp_{q_{i}^{-1}}(q_{i}^{-1}E_{i}q_{i}^{-H_{i}})
q_{i}^{H_{i}(H_{i}+2)/2}\\
=&
\exp_{q_{i}^{-1}}(-q_{i}^{-1}F_{i}q_{i}^{H_{i}})
\exp_{q_{i}^{-1}}(E_{i})
\exp_{q_{i}^{-1}}(-q_{i}^{-1}F_{i}q_{i}^{H_{i}})
q_{i}^{H_{i}(H_{i}+2)/2}
\end{split}
\end{equation}
where $E_{i},F_{i},H_{i}$ are the generators of $\Uhg$
corresponding to the simple root $\alpha_{i}$, $q_{i}=
q^{\<\alpha_{i},\alpha_{i}\>\hbar}$ and the $q$--exponential
is defined by
$$\exp_{q}(X)=
\sum_{n\geq 0}\frac{q^{n(n-1)/2}}{[n]_{q}!}X^{n}$$
where the $q$--factorials are given by
$$[n]_{q}=\frac{q^{n}-q^{-n}}{q-q^{-1}}
\quad\text{and}\quad
[n]_{q}!=[n]_{q}[n-1]_{q}\cdots [1]_{q}$$

Kwon investigated the $q$--Weyl group action of Artin's
braid group $B_{n}=B_{\sl{n}}$ on the zero weight space
of a simple $\Uhsl{n}$--module $\V$. He gave a general
criterion for it to be irreducible \cite{Kw}\footnote
{when coefficients are extended to the
field $\IC((\hbar))$ of formal Laurent series, which
we tacitly assume.} and showed moreover that this
criterion holds for all representations of $\Uhsl{3}$.
From these findings, he and Lusztig conjectured that
the action of $B_{n}$ on $\V[0]$ is irreducible for
any simple $\Uhsl{n}$--module $\V$. We shall prove
in this section the following 

\begin{theorem}\label{th:KL conj}
The Kwon--Lusztig conjecture is false for any complex,
simple Lie algebra $\g$ not isomorphic to $\sl{2},\sl
{3}$.
\end{theorem}

\subsection{Classical and Quantum Chevalley involution}

Let $\Th$ be the quantum Chevalley involution, \ie
the algebra automorphism of $\Uhg$ defined by
$$\Th(E_{i})=-F_{i},\quad
\Th(F_{i})=-E_{i}\quad\text{and}\quad
\Th(H_{i})=-H_{i}$$
As in the classical case, $\Th$ acts on any self--dual
\fd representation of $\Uhg$ leaving its zero weight
space invariant. Since $H_{i}$ acts as zero on $\V[0]$,
we see from \eqref{eq:triple qexp} that $\Th$ centralises
$\Bg$ on $\V[0]$. Corollary \ref{co:qTheta} below relates
the action of $\Th$ on $\V[0]$ to that of the classical
Chevalley involution $\Theta$ on $V[0]$. It will be used
to deduce theorem \ref{th:KL conj} from the results of
section \ref{se:Casimir}.

\begin{proposition}\label{pr:equiv}
There exists an algebra isomorphism $\Psi:\Uhg\rightarrow
\Ug\fml$ which is equal to the identity mod $\hbar$, acts as
the identity on $\h$ and satisfies
\begin{equation}\label{eq:equiv}
\Psi\circ\Th\circ\Psi^{-1}=\Theta
\end{equation}
\end{proposition}
\proof 
Let $\Phi:\Uhg\rightarrow\Ug\fml$ be an algebra isomorphism
equal to the identity mod $\hbar$ and acting as the identity on
$\h$ \cite [Prop. 4.3]{Dr2}. Set
$$\wt{\Theta}=\Phi\circ\Theta_\hbar\circ\Phi^{-1}=\Theta\circ\delta$$
where $\delta\in\id+\hbar\cdot\End(\Ug)\fml$ is such that $\delta
(p)=p$ for any $p\in U\h$. The involutivity of
$\wt{\Theta}$ yields $\Theta\circ\delta=\delta^{-1}\circ\Theta$
so that $a=\delta^{1/2}$ satisfies $a\circ\wt{\Theta}=\Theta\circ a$
and $\Psi=a\circ\Phi$ is the required isomorphism \halmos

\begin{corollary}\label{co:qTheta}
Let $\V$ be a self--dual, \fd $\Uhg$--module and let
$V=\V/\hbar\V$ be its reduction mod $\hbar$. Then,
$\V$ and $V\fml$ are isomorphic as $\h\rtimes\IZ_{2}
$--modules where the generator of $\IZ_{2}$ acts as
$\Th$ on $\V$ and as the classical Chevalley
involution $\Theta$ on $V$. In particular, $\Th$
acts as a scalar on $\V[0]$ iff it acts as a scalar
on $V[0]$.
\end{corollary}
\proof The isomorphism $\Psi:\Ug\fml\to\Ug\fml$
of proposition \ref{pr:equiv} endows $V\fml$ with
the structure of a $\Uhg$--module such that the
action of $\h\subset\Uhg$ coincides with that of
$\h\subset\g$. Since $\V$ and $V\fml$ have the
same reduction mod $\hbar$, they are isomorphic
as $\Uhg$, and therefore $\h$--modules. Equation
\eqref{eq:equiv} then guarantees that, under this
isomorphism, the Chevalley involution of $\V$ is
mapped to the Chevalley involution of $V$ \halmos\\

{\sc Proof of theorem \ref{th:KL conj}.} By proposition
\ref{pr:ext2g[0]} and theorem \ref{th:Casimir sln}, there
exists a simple, self--dual $\g$--module $V$ such that
$\Theta$ does not act as a scalar on $V[0]$. By corollary
\ref{co:qTheta}, $\Th$ does not act as scalar on $V\fml[0]$
and the latter is reducible under the $q$--Weyl group
action of $\Bg$ \halmos\\

We mention in passing the following $q$--analogue
of proposition \ref{pr:basic criterion}

\begin{proposition}\label{pr:qKostant}
Let $\V$ be a \fd representation of $\Uhg$ with
non--trivial zero weight space $\V[0]$. Then,
$\V[0]$ is irreducible under $\Bg$ iff it is
irreducible under $\Pg$.
\end{proposition}
\proof By proposition 1.2.1 of \cite{Sa},
$S_{i}$ acts on the zero--weight space of the
indecomposable $\Uhsl{2}^{\alpha_{i}}$--module
of dimension $2n+1$ as multiplication by $(-1)
^{n}q_{i}^{n(n+1)}$. It follows that the image
of $S_{i}$ in $\End(\V[0])$ is a polynomial in
the image of $S_{i}^{2}\in\Pg$ whence the
conclusion \halmos

\subsection{A Kohno--Drinfeld theorem for $q$--Weyl groups}

Let $\V$ be a \fd representation of $\Uhg$ and $V=\V/\hbar\V$
its reduction mod $\hbar$. It was conjectured in \cite{TL2}, by
analogy with the Kohno--Drinfeld theorem, that the $q$--Weyl
group action of $\Bg$ on $\V$ is equivalent to to the monodromy
action of $\Bg$ on $\nablak$ defined in the present paper\footnote
{this conjecture was independently formulated by De Concini
around 1995 (unpublished).}. This
conjecture is proved in \cite{TL3} for a number of pairs $(\g,V)
$ including vector representations of classical Lie algebras and
adjoint representations of all simple Lie algebras and in \cite
{TL2} for all representations of $\g=\sl{n}$. More precisely,

\begin{theorem}[\cite{TL2}]\label{th:monodromy}
Assume that $\g\cong\sl{n}$. Let $\mu$ be a weight of $V$
and
$$V^{\mu}=\bigoplus_{\nu\in W\mu}V[\nu]$$
the direct sum of the weight spaces of $V$ corresponding
to the Weyl group orbit of $\mu$. Let $\sigma(\Bg)\subset
N(T)$ be a Tits extension and
$$\rho^{\sigma}:\Bg\rightarrow GL(V^{\mu}\fmll)$$
the corresponding monodromy representation defined by
proposition \ref{th:existence} by regarding $h$ as a
formal variable. Let $\pi_{W}:\Bg\rightarrow GL(\V^
{\mu})$ be the $q$--Weyl group action. Then, $\rho_{h}$
and $\pi_{W}$ are equivalent for $\hbar=2\pi ih$.
\end{theorem}

Combining the above theorem with corollary \ref{co:gen red V[0]},
we obtain the following

\begin{proposition}
Let $\V$ be a \fd representation of $\Uhsl{n}$ and set
$V=\V/\hbar\V$. The following statements are equivalent
\begin{enumerate}
\item $\V[0]$ is (absolutely) irreducible under the $q$-Weyl group action of $B_n$.
\item $\V[0]$ is (absolutely) irreducible under the $q$-Weyl group action of $P_n$.
\item $V[0]$ is irreducible under the Casimir algebra $\Csl{n}$.
\end{enumerate}
\end{proposition}

In particular, by theorem \ref{th:Csl3}, the $q$--Weyl
group action of $P_{3}$ on the zero weight space $\V[0]$
of a $\Uhsl{3}$--module is always irreducible, a slight
refinement of a result of Kwon asserting the irreducibility
of $\V[0]$ under the full braid group $B_{3}$.

\section{Irreducible representations of $\Cg$}\label{se:irreps}

The aim of this section is to show that the connection
$\nablak$ yields irreducible monodromy representations
of $\Bg$ of arbitrarily large dimension. For $\g=\sl{n}$,
we show for example \S \ref{ss:irred A} that the weight
spaces of all Cartan powers of the adjoint representation
are irreducible under the Casimir algebra $\Cg$. For
$\g\ncong\sl{n}$, we obtain in \S \ref{ss:irred BCDEFG}
a somewhat weaker result : for every $p\in\IN$, the zero
weight space of the $p$th Cartan power of $\ad(\g)$ has
a subspace $K_{p}$ which is irreducible under $\Cg$
and such that $\lim_{p\rightarrow+\infty}\dim K_{p}=
+\infty$.

\subsection{Irreducible representations of $\C_{\sl{n}}$}\label{ss:irred A}

\begin{theorem}\label{th:sl Cartan}
For any $p,q\in\IN$, the action of $\Csln$ on the weight
spaces of the simple $\sl{n}$--module of highest weight
$(p,0,\ldots,0,-q)$ is irreducible.
\end{theorem}
\proof
For any $2\leq k\leq n$ and $a,b\in\IN$, set $\lambda^{(k)}_{a
,b}=(a,0,\ldots,0,-b)\in\IZ^{k}$ so that the \GZ basis of the
irreducible $\gl{n}$--module $V$ with highest weight $\lambda^{(n)}
_{p,q}$ is parametrised by patterns of the form
\begin{equation}\label{eq:pattern}
\lambda=
\begin{pmatrix}
\lambda^{(n)}_{p_{n},q_{n}}\\
\vdots\\
\lambda^{(k)}_{p_{k},q_{k}}\\
\vdots\\
\lambda^{(2)}_{p_{2},q_{2}}\\
r
\end{pmatrix}
\end{equation}
where the $p_{k}=p_{k}(\lambda)$, $q_{k}=q_{k}(\lambda)$ and
$r=r(\lambda)$ are integers satisfying
\begin{gather*}
p=p_{n}\geq p_{n-1}\geq\cdots\geq p_{2}\geq 0\\
q=q_{n}\geq q_{n-1}\geq\cdots\geq q_{2}\geq 0\\
p_{2}\geq r\geq -q_{2}
\end{gather*}

The vectors of a given weight $\mu=(\mu_{1},\ldots,\mu_{n})$
correspond to patterns satisfying in addition
$$r=\mu_{1}\quad\text{and, for any $2\leq k\leq n$,}\quad p_{k}-q_{k}=M_{k}$$
where $M_{k}=\sum_{i=1}^{k}\mu_{k}$. We claim that the
commuting Casimir operators $C_{\gl{k}}$, $2\leq k\leq
n$, have joint simple spectrum on the weight space $V
[\mu]$, with corresponding diagonal basis given by \GZ
vectors. Indeed, $C_{\gl{k}}$ acts on $v_{\lambda}\in
V[\mu]$, with $\lambda$ of the form \eqref{eq:pattern},
as multiplication by
\begin{equation}\label{eq:glk Cas}
\begin{split}
\<\lambda^{(k)}_{p_{k},q_{k}},
  \lambda^{(k)}_{p_{k},q_{k}}+2\rho^{(k)}\>
&=
p_{k}^{2}+q_{k}^{2}+(k-1)(p_{k}+q_{k})\\
&=
2q_{k}^{2}+2q_{k}(M_{k}+k-1)+M_{k}(M_{k}+k-1)
\end{split}
\end{equation}
where $2\rho^{(k)}=\sum_{i=1}^{k}\theta_{i}(k-2i+1)$ is the
sum of the positive roots of $\gl{k}$. Since $2q_{k}+M_{k}=
q_{k}+p_{k}\geq 0$ and the the right--hand side of \eqref
{eq:glk Cas} is a parabola with vertex at
$$q_{k}^{0}=-\half{1}(M_{k}+k-1)<-\half{M_{k}}$$
the $C_{\gl{k}}$--eigenvalue of a pattern $\lambda$ of form
\eqref{eq:pattern} and weight $\mu$ determines $p_{k}(\lambda)$
and $q_{k}(\lambda)$ uniquely as claimed.\\

We claim now that if $K\subseteq V[\mu]$ is a non--zero subspace
invariant under $\C_{\sl{n}}$, then $K=V[\mu]$. To see this, it
suffices to show that, for any given pattern $\lambda$ of the
form \eqref{eq:pattern} and weight $\mu$, there exists a \GZ
vector lying in $K$ such that the corresponding pattern has
the same $(n-1)$ row $\lambda^{(n-1)}_{p_{n-1},q_{n-1}}$ as
$\lambda$, for then a descending induction on $n$ shows that
$K$ contains all \GZ vectors of weight $\mu$. Since the 
Casimirs $C_{\gl{k}}$ have simple spectrum on $V[\mu]$,
$K$ contains at least one \GZ vector $v_{\lambda'}$. Let
$\lambda^{(n-1)}_{p'_{n-1},q'_{n-1}}$ be the $n-1$ row
of the corresponding pattern. If $p_{n-1}'=p_{n-1}$, then,
$$q_{n-1}'=
-M_{n-1}+p_{n-1}'=
-M_{n-1}+p_{n-1}=
q_{n-1}$$
and we are done. If $p_{n-1}'<p_{n-1}(\lambda)$, we may further
assume that $p_{n-1}'=\max_{\wt\lambda}p_{n-1}(\wt\lambda)$,
where the maximum is taken over all patterns $\wt\lambda$
such that $v_{\wt\lambda}\in K$ and $p_{n-1}(\wt\lambda)\leq
p_{n-1}$. Note then that
$$q_{n-1}'=
-\sum_{i=1}^{n-1}\mu_{i}+p_{n-1}'<
-\sum_{i=1}^{n-1}\mu_{i}+p_{n-1}=
q_{n-1}\leq q_{n}$$
It therefore follows by corollary \ref{co:nonzero} that
\begin{equation*}
\begin{split}
\kappa_{\theta_{n-1}-\theta_{n}}v_{\lambda'}
&=
(2E_{n-1,n}E_{n,n-1}+(E_{n-1,n-1}-E_{n,n}))v_{\lambda'}\\
&=
a v_{\lambda'+\varepsilon_{1}^{(n-1)}-\varepsilon_{n-1}^{(n-1)}}+
b v_{\lambda'-\varepsilon_{1}^{(n-1)}+\varepsilon_{n-1}^{(n-1)}}+
c v_{\lambda'}
\end{split}
\end{equation*}
for some $a,b,c\in\IC$ with $a\neq 0$. Hence, $K$ contains $v_
{\lambda'+\varepsilon_{1}^{(n-1)}-\varepsilon_{n-1}^{(n-1)}}$
in contradiction with the maximality of $p_{n-1}'$. The case 
$p_{n-1}'>p_{n-1}$ follows similarly \halmos\\

\remark A similar proof shows that the Casimir algebra $C_{\sl{n}}$
acts irreducibily on all weight spaces of the irreducible representations
with highest weight of the form $(p,q,0,\ldots,0)$ where $p,q\in\IZ$
satisfy $ p\geq q\geq 0$. More generally, one can
show that if the commuting  Casimirs $C_{\gl{k}}$, $k=2\ldots n$,
have joint simple spectrum on the weight space $V[\mu]$ of a
simple $\sl{n}$--module $V$, then $\C_{\sl{n}}$ acts irreducibly
on $V[\mu]$. The proof is similar to that of theorem \ref
{th:sl Cartan} but somewhat more involved technically and
will be given in a future publication.

\subsection{Irreducible representations of $\Cg$, $\g\ncong\sl{n}$}
\label{ss:irred BCDEFG}

Let $\g$ be a complex, simple Lie algebra not isomorphic to $\sl{n}$
and let $\theta$ be the highest root of $\g$.

\begin{theorem}
For any $p\in\IN$, there exists a subspace $K_{p}$ of the
zero weight space of the simple $\g$--module with highest
weight $p\theta$ which is irreducible under $\Cg$ and such
that $\lim_{p\rightarrow\infty}\dim K_{p}=\infty$.
\end{theorem}
\proof We shall need the following simple

\begin{lemma}
Let $R$ be a root system and $\alpha\neq\pm\beta\in R$ two
long roots which are not orthogonal. Then, $R\cap(\IZ\alpha
+\IZ\beta)$ is a root system of type $A_{2}$.
\end{lemma}
\proof Since $\|\alpha\|=\|\beta\|$, one has $\<\alpha,\beta
^{\vee}\>=\pm 1$. Replacing $\beta$ by $-\beta$ if necessary,
we may assume that $\<\alpha,\beta^{\vee}\>=-1$. Thus, $\pm
\alpha,\pm\beta,\pm(\alpha+\beta)=\pm\sigma_{\beta}\alpha\in
R$ and these are the only $\IZ$--linear combinations of $\alpha,
\beta$ which lie in $R$ since any other has norm strictly larger
than $\|\alpha\|$ \halmos\\

Assume first that $\g\ncong\sp{2n}$, $n\geq 2$. An inspection
of the tables in \cite{Bo} then shows that there is a unique simple
root $\alpha$ of $\g$ which is not orthogonal to $\theta$ and is,
moreover, long. Applying the above lemma to the pair $(\alpha,
\theta)$ yields a subalgebra
$$\ll=
\IC h_{\alpha}\oplus\IC h_{\theta}
\bigoplus_{\gamma\in R(\g)\cap(\IZ\alpha+\IZ\theta),\gamma\succ 0}
\IC e_{\gamma}\oplus \IC f_{\gamma} 
\subset\g$$
which is isomorphic to $\sl{3}$ and has as highest root vector
$e_{\theta}$. Choose $\h_{\ll}=\IC h_{\alpha}\oplus\IC h_{\theta}$
as Cartan subalgebra of $\ll$ and denote by $\ll^{*p}$ (resp.
$\g^{*p}$) the irreducible representation of $\ll$ (resp. $\g$)
with highest weight $p\theta$. Since $\ll^{*p}$ is generated by
$e_{\theta}^{\otimes p}$ inside $\ll^{\otimes p}\subset\g^{\otimes p}$,
it follows that $\g^{*p}$ contains $\ll^{*p}$ as $\ll$--submodule.
This inclusion induces one of weight spaces $\ll^{*p}[0]\subset
\g^{*p}[0]$ since $\h=\h_{\ll}\oplus\h_{\ll}^{\perp}$ and
$\h_{\ll}^{\perp}$ centralises $\ll$. By theorem \ref{th:Csl3},
$\ll^{*p}[0]$ is irreducible under $\C_{\ll}$. Let $U=\Cg\ll^{*p}[0]$
be the $\Cg$--submodule of $\g^{*p}[0]$ generated by $\ll^{*p}[0]$
and decompose it as a sum $\bigoplus_{i}U_{i}$ of irreducible
summands with projections $p_{i}$. By Schur's lemma, the
restriction of each $p_{i}$ to $\ll^{*p}[0]$ is either
zero or injective. Thus, the dimension of at least one of
the $U_{i}$'s is greater or equal to that of $\ll^{*p}[0]$
and therefore tends to infinity with $p$. If $\g\cong\sp{2n}$,
the \KT branching rules \cite[thm. A1]{KT} show that $\res_{\g}
^{\sl{n}}\g^{*p}\supset{\sl{n}}^{*p}$ and we may conclude as
above with $\l=\sl{n}$ \halmos

\section{Zero weight spaces of self--dual $\g$--modules}
\label{se:selfdual V[0]}

The aim of this section is to show that, when $\g$ is isomorphic
to $\sl{n}$, $n\geq 4$, or $\g_2$, the zero weight space of most
{\it self--dual}, simple $\g$--modules is reducible under the Casimir
algebra $\Cg$, thus strengthening the results of section \ref{se:Casimir}.
We classify the self--dual $V$ with $V[0]$ irreducible under $\C_{\sl{n}}$
and $\C_{\g_2}$ in \ref{ss:V[0] A} and \ref{ss:V[0] G} respectively.
Our calculations rely on the use of the Chevalley involution
$\Theta$ via proposition \ref{pr:Chevalley}. They show in fact
that $V[0]$ is irreducible under $\Cg$ if, and only if $\Theta$
acts as a scalar on it. We conjecture in \ref{ss:Chevalley conjecture}
that this holds for any simple Lie algebra $\g$.\\

\remark By corollary \ref{co:qTheta}, our results imply that the
$q$--Weyl group action of $\Bg$ on the zero weight spaces of
most self--dual $\Uhg$--modules is reducible.

\subsection{Self--dual representations of $\sl{n}$}
\label{ss:V[0] A}

\begin{theorem}\label{th:V[0] A}
Let $V$ be a simple, self--dual $\sl{n}$--module
with non--trivial zero weight space $V[0]$. Then,
$V[0]$ is irreducible under the Casimir algebra
$\Csl{n}$ if the highest weight $\lambda$ of $V$
is of one of the following forms
\begin{enumerate}
\item $\lambda=(p,0,\ldots,0,-p)$, $p\in\IN$.
\item $\lambda=(\underbrace{1,\ldots,1}_{k},
0,\ldots,0, \underbrace{-1,\ldots,-1}_{k})$,
$0\leq k\leq n/2$
\item $\lambda=(p,p,-p,-p)$, $p\in\IN$.
\end{enumerate}
Conversely, if $\lambda$ is of none of the above
forms, there exists a $k<n$ and a simple $\gl{k}
$--summand $U\subset V$ with $U\ncong U^*$
and $U\cap V[0]\neq\{0\}$. In particular, $V[0]$ is
reducible under $\Csl{n}$.
\end{theorem}

Since the case (i) follows from theorem \ref{th:sl Cartan},
the ``if'' part of theorem \ref{th:V[0] A} is settled by
the following two lemmas.

\begin{lemma}\label{le:Vnk}
Let $V_{n,k}$ be the simple $\sl{n}$--module with highest
weight
$$\lambda\nn_{k}=
(\underbrace{1,\ldots,1}_{k},0,\ldots,0,
 \underbrace{-1,\ldots,-1}_{k})$$
Then $V_{n,k}[0]$ is irreducible under $\C_{\sl{n}}$.
\end{lemma}
\proof We claim that the Casimirs $C_{\gl{m}}$, $m=2,
\ldots n$, have joint simple spectrum on $V_{n,k}[0]$.
Indeed, the $m$th row of a zero weight \GZ pattern
corresponding to $V_{n,k}$ is of the form $\lambda^
{(m)}_{l}$, for some $0\leq l\leq m/2$. Since $C_
{\gl{m}}$ acts as multiplication by $2l(m-l+1)$
on the representation with highest weight $\lambda
^{(m)}_{l}$ and the function $f(x)=2x(m-x+1)$ is
injective on the interval $[0,(m+1)/2]$, the $C_{\gl{m}}
$--eigenvalue of a zero weight \GZ pattern determines
its $m$th row uniquely, as claimed. The proof is now
completed as in theorem \ref{th:sl Cartan} \halmos

\begin{lemma}\label{le:C2pExt2}
For any $p\in\IN$, let $V_{p}$ be the simple $\sl{4}
$--module with highest weight $(p,p,-p,-p)$. Then,
$V_{p}[0]$ is irreducible under $\Csl{4}$.
\end{lemma}
\proof The \GZ patterns corresponding to $V_{p}[0]$ are
of the form
$$\begin{array}{ccccccc}
p& &p& &-p& &-p\\
 &p& &0& &-p&  \\
 & &q& &-q& &  \\
 & & &0& & &   \\
\end{array}$$
for some $0\leq q\leq p$, so that they are separated by
the action of the Casimir of $\sl{2}\subset\sl{n}$. The
proof is now completed as in theorem \ref{th:sl Cartan}
\halmos\\

{\sc Proof of theorem \ref{th:V[0] A}.} Assume that $\lambda$
is not of the form (i)--(iii) and let
$$s(\lambda)=|\{i=1\ldots n-1|\lambda_{i}-\lambda_{i+1}>0\}|$$
be the number of steps in the corresponding Young diagram.
Noting that $s(\lambda)=2$ for the highest weights of the
form (i)--(iii), we begin by proving that $V[0]$ is reducible
under $\Csl{n}$ if $s(\lambda)\geq 3$. Suppose first that
$n=2k+1$ is odd. Since $V[0]\neq\{0\}$, we may assume that
the sum $|\lambda|$ of the entries in $\lambda$ is zero so
that, by self--duality, $\lambda$ is of the form
$$\lambda=(a_1,\ldots,a_l,0,\ldots,0,-a_l,\ldots,-a_1)$$
for some $a_{1}\geq\cdots\geq a_{l}>0$, with at least one
middle zero. Since $s(\lambda)\geq 3$, there exists some
$1\leq i\leq l-1$ such that $a_{i}>a_{i+1}$. Let $\mu$
be the $\sl{n-1}$--weight obtained by replacing $a_{i}$
by $a_{i}-1$ in the $i$th position, $-a_{l}$ by $-a_{l}+
1$ in the $n-l+1$th position and by omitting the middle
zero. Then, $|\mu|=\sum_{i=1}^{n-1}\mu_{i}=|\lambda|=0$
and
$$\mu_{i}+\mu_{n-i}=-1\neq 1=\mu_{l}+\mu_{n-l}$$
so that $\mu$ is a non--self dual weight of $\sl{n-1}$
such that $V_{\mu}\cap V[0]\neq\{0\}$.\\

Consider now the case where $n=2k$ is even. Assuming
again that $|\lambda|=0$, we find that
$$\lambda=
(a_1,\ldots,a_l,\underbrace{0,\ldots,0}_{n_{0}},-a_l,\ldots,-a_1)$$
where the number $n_{0}$ of zeroes is even. If $n_{0}>0$,
the $\sl{n-1}$ weight $\mu$ obtained by replacing the two
middle zeroes by a single one in $\lambda$ is self--dual
and satisfies $|\mu|=0$ and $s(\mu)=s(\lambda)$. By our
previous analysis, there therefore exists a non--self dual
$\sl{n-2}$ weight $\nu\prec\mu\prec\lambda$ such that
$|\nu|=0$ and $V_{\nu}\cap V[0]\neq 0$. If, on the other
hand, $n_{0}=0$, then
$$\lambda=(a_1,\cdots,a_{n/2},-a_{n/2},\cdots,-a_1)$$
with $a_{n/2}>0$. Let $1\leq i\leq n/2-1$ be such that
$a_{i}>a_{i+1}$. Then, the non--self dual $\sl{n-1}$--weight
$\mu\prec\lambda$ obtained by replacing $a_{i}$ by $a_{i}-1$
in the $i$th position and the pair $a_{n/2},-a_{n/2}$ allows
to conclude.\\

Consider now the case $s(\lambda)\leq 2$. By self--duality,
we may take $\lambda$ of the form
$$\lambda=
(\underbrace{p,\ldots,p}_{k},0,\ldots,0,
 \underbrace{-p,\ldots,-p}_{k})$$
for some $0\leq k\leq n/2$. By assumption, $k\geq 2$ and
$p>1$ since $\lambda$ is not of the forms (i)--(ii). If
$n$ is odd, replacing the innermost pair $(p,-p)$ 
by $(p-1,-(p-1))$ and suppressing the middle zero yields
an $\sl{n-1}$ weight $\mu$ with $s(\mu)\geq 3$ and our
previous analysis allows to conclude. If $n$ is even and
there are two or more middle zeroes, we suppress one of
them to obtain an $\sl{n-1}$--weight $\nu$ of of the form
treated in the previous paragraph. If there are no middle
zeroes, then by assumption $n\geq 6$. We change the innermost
pair $(p,-p)$ to $0$ and proceed as above \halmos

\subsection{Representations of $\g_{2}$}
\label{ss:V[0] G}

Recall that, for $\g=\g_{2}$, every $\g$--module $V$
is self--dual and has a non--trival zero weight space.
The aim of this subsection is to prove the following.

\begin{theorem}\label{th:V[0] G}
The zero weight space $V[0]$ of a simple $\g_{2}
$--module $V$ is irreducible under the Casimir
algebra $\C_{\g_{2}}$ iff $V$ is a trivial or
fundamental representation, or its second Cartan
power.
\end{theorem}

The proof of the theorem is given in the next three
propositions. We begin by reviewing Perroud's branching
rules for the equal rank inclusion $\sl{3}\subset\g_{2}
$ \cite{Pe}. Let $\alpha_{1},\alpha_{2}$ be the long and
short simple roots of $\g_{2}$ respectively and $\varpi_
{1},\varpi_{2}$ the corresponding fundamental weights\footnote
{we follow Perroud's convention \cite{Pe} which are
the opposite of the usual ones \cite{Bo,FH}}. Let $V_
{\lambda}$ be the simple $\g_{2}$--module with highest
weight $\lambda=m_{1}\varpi_{1}+m_{2}\varpi_{2}$. Consider
the set of \GZ patterns $\mu(a,b,c)$ of the form
\begin{equation}\label{eq:P patterns}
\begin{array}{ccccc}
m_{1}+m_{2}& &m_{2}& &0\\
	   &a&	   &b& \\
	   & &c	   & &
\end{array}
\end{equation}
Then,
\begin{equation}\label{eq:Perroud}
\res_{\g_{2}}^{\sl{3}}V_{\lambda}=
\bigoplus_{\mu(a,b,c)}
V_{(m_{1}+c,a-m_{2}+b,0)}
\end{equation}

\begin{lemma}\label{le:ad[0]}
For any complex, simple Lie algebra $\g$, $\h=\ad(\g)
[0]$ is irreducible under $\Cg$.
\end{lemma}
\proof This follows from proposition \ref{pr:basic criterion}
since $\h$ is irreducible under $W$ \halmos\\

\begin{proposition}
The zero weight spaces of the fundamental representations
of $\g$ and of their second Cartan powers are irreducible
under $\C_\g$.
\end{proposition}
\proof The zero weight space of $V_{\varpi_2}\cong\IC^7$ is
one--dimensional and therefore irreducible under $\Cg$. The
branching rules \eqref{eq:Perroud} yield
$$\res_{\g_{2}}^{\sl{3}} V_{2\varpi_2}=\ad(\sl{3})\oplus R$$
where $R$ is a reducible $\sl{3}$--modules with trivial zero
weight space. The irreducibility of $V_{2\varpi_2}[0]$ therefore
follows from lemma \ref{le:ad[0]} for $\g=\sl{3}$ or from theorem
\ref{th:Csl3}. Since $V_{\varpi_1}\cong\ad(\g_2)$, the irreducibility
of $V_{\varpi_1}[0]$ under $\Cg$ follows from lemma \ref{le:ad[0]}.
Finally, by \eqref{eq:Perroud},
$$\res_\g^{\sl{3}}C^2\ad(\g)=C^2\ad(\sl{3})\oplus V\oplus R$$
where
\begin{align*}
V&=V_{(2,1,0)}\cong\ad(\sl{3})\\
R&=
V_{(3,2,0)}\oplus V_{(3,1,0)}\oplus
V_{(2,2,0)}\oplus V_{(2,0,0)}
\end{align*}
and $R$ has a trivial 0 weight space. One readily checks that $C^{2}\ad
(\sl{3})$ and $V$ are distinguished by the Casimir eigenvalue of $\sl{3}$
so that, by theorem \ref{th:Csl3}, the zero weight spaces of $C^{2}\ad(\sl
{3})$ and $V$ are irreducible and inequivalent representations of $\C_{\sl
{3}}$. Since $C^{2}\ad(\g)[0]=C^{2}\ad(\sl{3})[0]\oplus V[0]$ it suffices
to show that $C^{2}\ad(\sl{3})[0]$ is not invariant under $\Cg$. Let $\theta$
be the highest root of $\g$ and $e_{\theta},f_{\theta},h_{\theta}$ a corresponding
$\sl{2}^{\theta}$ triple. This triple lies in $\sl{3}$ so that
$$v_{\theta}=
-1/4\ad(f_{\theta})\ad(e_{\theta}) e_{\theta}^{2}=
e_{\theta}\cdot f_{\theta}-\half{1}h_{\theta}^{2}
\in C^{2}\ad(\sl{3})\subset C^{2}\ad(\g)$$
where we are realising $C^{2}\ad(\g)$ as the highest
weight component of $S^{2}\g$. Let $\alpha$ be a short root of
$\g$ such that $\<\theta^{\vee},\alpha\>=1$ so that the $\alpha
$--string through $\theta$ is of the form $\theta-2\alpha,\theta
-\alpha,\theta$. It is easy to see that such an $\alpha$
exists by consulting the tables in \cite{Bo}. A simple computation
using a Chevalley basis of $\g$ yields
\begin{align*}
-\half{1}\ad(f_{\alpha})\ad(e_{\alpha})\medspace
h_{\theta}^{2}&=
e_{\alpha}\cdot f_{\alpha}-h_{\alpha}\cdot h_{\theta}\\
\ad(f_{\alpha})\ad(e_{\alpha})\medspace
e_{\theta}\cdot f_{\theta}&=
\pm e_{\theta-\alpha}\cdot f_{\theta-\alpha}
\pm 2 e_{\theta}\cdot f_{\theta}
\end{align*}

where the signs depend on the choice of the root vectors.
Thus,
$$\frac{1}{\<\alpha,\alpha\>}C_{\alpha} v_{\theta}=
\ad(f_{\alpha})\ad(e_{\alpha})\medspace v_{\theta}=
e_{\alpha}\cdot f_{\alpha}\pm
e_{\theta-\alpha}\cdot f_{\theta-\alpha}\pm
2 e_{\theta}\cdot f_{\theta}-
h_{\alpha}\cdot h_{\theta}$$
which does not lie in $C^{2}\ad(\sl{3})$ since $\alpha$
is short \halmos

\begin{proposition}
Let $V$ be a simple $\g_{2}$--module with highest
weight $\lambda=m_{1}\varpi_{1}+m_{2}\varpi_{2}$.
If $m_{1}+m_{2}\geq 3$, the zero weight space of
$V$ is reducible under $\C_{\g_{2}}$.
\end{proposition}
\proof By proposition \ref{pr:Chevalley}, it suffices
to prove that the restriction of $V$ to $\sl{3}$ 
contains an irreducible summand $U$ with $U[0]\neq
\{0\}$ and $U\ncong U^{*}$. We begin by treating
the special cases $m_{1}=0$ and $m_{2}=0$. Assume
first that $m_{2}=0$ so that $b=0$ in \eqref{eq:P patterns}.
If $m_{1}=0 \mod 3$, then setting $a=c=0$ in \eqref
{eq:Perroud} yields $\res_{\g_{2}}^{\sl{3}}V\supset
V_{(m_{1},0,0)}$. Similarly, if $m_{1}=1\mod 3$,
with $m_{1}>1$, taking $a=1$ and $c=1$ yields
$\res_{\g_{2}}^{\sl{3}} \supset V_{(m_{1}+1,1,0)}$.
Finally, if $m_{1}=2\mod 3$, $m_{1}>2$, $a=1,c=0$
yields $\res_{\g_{2}}^{\sl{3}}\supset V_{(m_{1},1,0)}$
as required. Assume now that $m_{1}=0$ and $m_{2}
\geq 3$ so that $a=m_{2}$ in \eqref{eq:P patterns}.
Then, taking $b=0,c=3$, we find $\res_{\g_{2}}^{\sl{3}}
\supset V_{(3,0,0)}$.\\

Consider now the case $m_{1},m_{2}>0$. The values
of $(a,b,c)$ corresponding to the \GZ patterns \eqref
{eq:P patterns} are readily seen to span the integral
points of a convex polytope in $\IR^{3}$ with vertices
given by
\begin{gather*}
(m_{2},0,0),(m_{2},0,m_{2}),(m_{2},m_{2},m_{2}),\\
(m_{1}+m_{2},0,0),(m_{1}+m_{2},0,m_{1}+m_{2}),\\
(m_{1}+m_{2},m_{2},m_{2}),(m_{1}+m_{2},m_{2},m_{1}+m_{2})
\end{gather*}
The image $P(m_{1},m_{2})\subset\IR^{3}$ of this polytope
under the Perroud map $\pi:(a,b,c)\rightarrow (m_{1}+c,a-
m_{2}+b,0)$ is the convex hull of the images of the above
points, namely
\begin{gather*}
(m_{1},0,0),(m_{1}+m_{2},0,0),(m_{1}+m_{2},m_{2},0),\\
(m_{1},m_{1},0),(2m_{1}+m_{2},m_{1},0),\\
(m_{1}+m_{2},m_{1}+m_{2},0),(2m_{1}+m_{2},m_{1}+m_{2},0)
\end{gather*}
and is readily seen to be described by the following
inequalities in the plane $(\mu_{1},\mu_{2},0)\subset
\IR^{3}$
\begin{gather}
m_{1}\leq\mu_{1}\leq 2m_{1}+m_{2}
\label{eq:ineq 1}\\
0\leq\mu_{2}\leq m_{1}+m_{2}
\label{eq:ineq 2}\\
0\leq\mu_{1}-\mu_{2}\leq m_{1}+m_{2}
\label{eq:ineq 3}
\end{gather}
Moreover,
$$\res_{\g_{2}}^{\sl{3}} V=
\bigoplus_{\mu\in P(m_{1},m_{2})\cap\IN^{3}}
V_{\mu}\otimes\IC^{|\pi^{-1}(\mu)|}$$
We seek to derive a contradiction from the assumption
that all summands with non--trivial zero weight space are
self--dual. Let $U$ be a summand with $U[0]\neq\{0\}$
and $U\cong U^{*}$ so that its highest weight is of the
form $\mu=(2k,k,0)$ for some $k\in\IN$. Let $\tau:\IZ^
{3}\rightarrow\IZ^{3}$ be defined by
$$\tau(\nu_{1},\nu_{2},\nu_{3})=
(\nu_{1}-1,\nu_{2}+1,\nu_{3})$$
so that, if $\tau(\mu)\in P(m_{1},m_{2})$ (resp. $\tau^{-1}
(\mu)\in P(m_{1},m_{2})$) then $V_{\tau(\mu)}\subset V$
(resp. $V_{\tau^{-1}(\mu)}\subset V$) is a non self--dual
$\sl{3}$--summand with non--trivial zero weight space. We
shall need the following

\begin{lemma}
Assume that $m_{1},m_{2}\neq 0$ and that $\mu=(2k,k,0)\in
P(m_{1},m_{2})$.
\begin{enumerate}
\item[(i)] If $\tau(\mu)\notin P(m_{1},m_{2})$, then
$k\in\{1,m_{1}/2\}$.
\item[(ii)] If $\tau^{-1}(\mu)\notin P(m_{1},m_{2})$,
then $k\in\{m_{1}+m_{2}-1,m_{1}+m_{2}/2\}$.
\end{enumerate}
\end{lemma}
\proof (i) By assumption, $\tau(\mu)$ violates at least
one of the inequalities \eqref{eq:ineq 1}--\eqref{eq:ineq 3},
so that at least one of the following conditions holds
\begin{gather*}
\mu_{1}=m_{1}\\
\mu_{2}=m_{1}+m_{2}\\
\mu_{1}-\mu_{2}\in\{0,1\}
\end{gather*}
The condition $k=\mu_{1}-\mu_{2}=0$ is ruled out by
the fact that $(0,0,0)\notin P(m_{1},m_{2})$ if $m_{1}>0$.
Similarly, $k=\mu_{2}=m_{1}+m_{2}$ leads to $(2(m_{1}+
m_{2}),m_{1}+m_{2},0)\in P(m_{1},m_{2})$ which violates
\eqref{eq:ineq 1} since $m_{2}>0$. We are therefore
left with $\mu_{1}=m_{1}$ or $\mu_{1}-\mu_{2}=1$ which
lead to $k=m_{1}/2,1$ respectively. (ii) Similarly,
$\tau^{-1}(m)\notin P(m_{1},m_{2})$ iff at least one
of the following equations holds
\begin{gather*}
\mu_{1}=2m_{1}+m_{2}\\
\mu_{2}=0\\
\mu_{1}-\mu_{2}\in\{m_{1}+m_{2}-1,m_{1}+m_{2}\}
\end{gather*}
$\mu_{2}=0$ and $\mu_{1}-\mu_{2}=m_{1}+m_{2}$ imply
that $\mu=(0,0,0)$ and $\mu=(2(m_{1}+m_{2}),m_{1}+
m_{2})$ respectively both of which are ruled out by
$m_{1},m_{2}>0$. Thus, $\mu_{1}=2m_{1}+m_{2}$ or
$\mu_{1}-\mu_{2}=m_{1}+m_{2}-1$ hold yielding
$k\in\{m_{1}+m_{2}/2,m_{1}+m_{2}-1\}$ \halmos\\

Returning to our main argument, if all $\sl{3}$--summands
in $V$ with non--trivial zero weight spaces are self--dual
then $\tau(\mu),\tau^{-1}(\mu)\notin P(m_{1},m_{2})$ for
any $\mu\in P(m_{1},m_{2})$ of the form $(2k,k,0)$. By
the above lemma, this implies
$$
\{1,m_{1}/2\}\cap\{m_{1}+m_{2}-1,m_{1}+m_{2}/2\}\neq\emptyset
$$
so that at least one of the following equations holds
\begin{gather*}
m_{1}+m_{2}=2\\
m_{1}+m_{2}/2=1\\
m_{1}/2+m_{2}=1\\
m_{1}/2+m_{2}/2=0
\end{gather*}
contradicting the fact that $m_{1}+m_{2}\geq 3$ \halmos

\begin{proposition}\label{pr:g211}
Let $U$ be the Cartan product of the two fundamental
representations of $\g_{2}$. Then, $U[0]$ is reducible
under $\C_{\g_{2}}$.
\end{proposition}
\proof It suffices to show that the Chevalley
involution of $\g=\g_{2}$ does not act as a
scalar on $U[0]$. Let $\alpha_{1},\alpha_{2}$
be the short and long simple roots respectively
\footnote{we adhere now to the standard notation
\cite{FH}} and label the positive roots by
$$
\alpha_{i}=(i-2)\alpha_{1}+\alpha_{2},
\thickspace\thickspace 2\leq i\leq 5
\qquad\text{and}\qquad
\alpha_{6}=3\alpha_{1}+2\alpha_{2}
$$
so that the highest root is $\theta=\alpha_{6}$.
The corresponding fundamental weights $\varpi_{1},
\varpi_{2}$ of $\g$ are
$$
\varpi_{1}=2\alpha_{1}+\alpha_{2}=\alpha_{4}
\qquad\text{and}\qquad
\varpi_{2}=3\alpha_{1}+2\alpha_{2}=\alpha_{6}
$$
so that $V_{\varpi_{2}}\cong\ad(\g_{2})$ and
$V=V_{\varpi_{1}}\cong\IC^{7}$ has weights $
0$ and $\pm\alpha_{i}$, $i=1,3,4$ \cite[\S
22.1]{FH}. Choose a Cartan--Weyl basis $e_{
\alpha_{i}},f_{\alpha_{i}},h_{\alpha_{1}},
h_{\alpha_{2}}$ of $\g$ and a weight basis
$v_{\pm\alpha_{i}}$, $i=1,3,4$ and $v_{0}$
of $V_{\varpi_{1}}$ where $v_{\beta}$ has
weight $\beta$. The highest weight vector
in $U\subset\g_{2}\otimes V$ is $e_{\alpha
_{6}}\otimes v_{\alpha_{4}}$ so that
$$
u_{0}=
f_{\alpha_{4}}f_{\alpha_{6}}\medspace 
e_{\alpha_{6}}\otimes v_{\alpha_{4}}\in U[0]
$$

Computing $u_{0}$ explicitly yields 
\begin{equation*}
\begin{split}
u_{0}
&=
 f_{\alpha_{4}}\left(
-h_{\alpha_{6}}\otimes v_{\alpha_{4}}+
ae_{\alpha_{6}}\otimes v_{-\alpha_{3}}\right)\\
&=
-\alpha_{4}(h_{\alpha_{6}})f_{\alpha_{4}}\otimes v_{\alpha_{4}}
+b h_{\alpha_{6}}\otimes v_{0}
+c e_{\alpha_{3}}\otimes v_{-\alpha_{3}}
\end{split}
\end{equation*}

where the constants $a,b,c$ depend upon the choices
of the basis of $V_{\varpi_{i}}$, $i=1,2$ and are
not zero by elementary $sl_{2}$--representation
theory. On the other hand, if $\Theta$ is the Chevalley
involution of $\g$, then $\Theta \medspace h_{\alpha_{i}}
=-h_{\alpha_{i}}$ and, up to multiplicative constants
\begin{xalignat*}{3}
\Theta \medspace e_{\alpha_{i}} &= f_{\alpha_{i}}&
\Theta \medspace f_{\alpha_{i}} &= e_{\alpha_{i}}&
\Theta \medspace v_{\alpha_{i}} &= v_{-\alpha_{i}}
\end{xalignat*}
so that $\Theta u_{0}$ is not proportional to $u_{0}$
\halmos

\subsection{Some conjectures}\label{ss:Chevalley conjecture}

Let us record the following corollary of the proofs of theorems
\ref{th:V[0] A} and \ref{th:V[0] G}.

\begin{theorem}
Let $\g$ be $\sl{n}$ or $\g_{2}$ and let $V$ be a simple,
self--dual $\g$--module with $V[0]\neq\{0\}$. Then, if
$V[0]$ is reducible under the Casimir algebra of $\g$,
the Chevalley involution does not act as a scalar on
$V[0]$.
\end{theorem}

In other words, for the above representations, the failure of
the Chevalley involution to act as a scalar on $V[0]$ is the
only mechanism which causes $V[0]$ to be reducible under
the Casimir algebra $\Cg$. It is therefore natural to make
the following

\begin{conjecture}
Let $\g$ be a complex, simple Lie algebra and let $V$
be a simple $\g$--module which is self--dual and has
a non--trivial zero weight space $V[0]$. Then $V[0]$
is irreducible under the Casimir algebra of $\g$ iff
the Chevalley involution of $\g$ acts as a scalar on
$V[0]$.
\end{conjecture}

We plan to address this conjecture in a future publication.
We note here that it holds for the irreducible representations
of $\g=\so{2n},\so{2n+1},\sp{2n}$ whose highest weight
$\lambda=(\lambda_1,\ldots,\lambda_n)$ satisfies $\lambda
_i=0$ for $i>n/2$ \cite{HMTL}.\\

It would also be desirable to formulate in a way independent
of the Lie type of $\g$ our observation that the zero weight
space of 'most' self--dual $\g$--modules is reducible under
the Casimir algebra. At the very least, for example, we make
the following

\begin{conjecture}
Let $\g\ncong\sl{2},\sl{3}$ and let $V$ be a simple, self--dual
$\g$--module with $V[0]\neq 0$. If the highest weight of $V$
is regular, then $V[0]$ is reducible under $\Cg$.
\end{conjecture}

The above conjecture is true for $\g=\sl{n}$ and $\g=\g_{2}$
by theorems \ref{th:V[0] A} and \ref{th:V[0] G}. It also holds
for $\g=\so{2n},\so{2n+1},\sp{2n}$ \cite{HMTL}.

\section{Appendix : The centraliser of the Casimir algebra}
\label{se:etingof}

The results in this section are due to P. Etingof \cite{Et}
to whom we are grateful for allowing us to reproduce
them here. Our aim is to prove the following

\begin{theorem}\label{th:centraliser}
The centraliser of the Casimir algebra $\Cg$ of $\g$ is
generated by the Cartan subalgebra $\h$ and the centre
$Z(\Ug)$ of $\Ug$.
\end{theorem}

An immediate corollary of the above result is the following
special case of a theorem of Knop \cite[thm. 10.1]{Kn}
\footnote{In an earlier version of this paper, Knop's result
was used, in conjuntion with theorem \ref{th:Zariski}, to
prove theorem \ref{th:centraliser}.}.

\begin{corollary}[Knop]\label{co:Knop}
The centre of $\Ugh$ is isomorphic to $Z(\Ug)\otimes U\h$.
\end{corollary}

\subsection{} The proof of theorem \ref{th:centraliser} will be
given in \S \ref{ss:centraliser}. It rests on the following result
which is of independent interest

\begin{theorem}\label{th:Zariski}
For any non--negative linear combination $\beta\in\bigoplus
\alpha_i\cdot\IN$ of simple roots, there exists a Zariski open
set $O_{\beta}\subset\h^{*}$ such that, for any $\mu\in O_{
\beta}$, the $(\mu-\beta)$--weight space $M_{\mu}[\mu-\beta]
$ of the Verma module with highest weight $\mu$ is irreducible
under the action of the Casimir algebra $\Cg$.
\end{theorem}
\proof We need a preliminary result. Fix a weight $\lambda\in
\h^{*}$ and let $M_{t^{2}\lambda}$ be the Verma module with
highest weight $t^{2}\lambda$, where $t\in\IC^{*}$ is some
non--zero complex number. Consider the standard
identifications

$$M_{t^{2}\lambda}
\stackrel{\imath}{\longrightarrow}
\Un_{-}\stackrel{\sigma^{-1}}{\longrightarrow}
\Sn_{-}$$

where $\sigma$ is the symmetrisation map. The corresponding
isomorphism $M_{t^{2}\lambda}\cong \Sn_{-}$ is one of
$\h$--modules provided the adjoint action of $\h$ on
$\Sn_{-}$ is tensored by the character $t^{2}\lambda$.
Denoting the generators of $\Sn_{-}$ by $x_{\alpha}$
and transporting the action of $\g$ on $M_{t^{2}\lambda}$
to $\Sn_{-}$, we have the following

\begin{lemma}\label{le:asymptotic}
Let $d$ be the grading operator on $\Sn_{-}$. Then, for any
$t\in\IC^{*}$
\begin{align}
t^{d}e_{\alpha}t^{-d}
&=t\cdot \<\lambda,\alpha^{\vee}\>\partial_{\alpha}+O(1)
\label{eq:e}\\
t^{d}f_{\alpha}t^{-d}
&=t\cdot x_{\alpha}+O(1)
\label{eq:f}\\
t^{d}h_{\alpha}t^{-d}
&=t^{2}\cdot\<\lambda,\alpha^{\vee}\>+O(1)
\label{eq:h}
\end{align}

where the terms $O(1)$ have a finite limit for $t
\longrightarrow\infty$.
\end{lemma}
\proof \eqref{eq:e} Let $v_{t^{2}\lambda}\in M_{t^{2}\lambda}$
be the highest weight vector. Then, for any sequence of positive
roots $\beta_{1},\ldots,\beta_{k}$, we have
\begin{multline*}
t^{d}e_{\alpha}t^{-d}
\medspace
x_{\beta_{1}}\cdots x_{\beta_{k}}
=\\
t^{d-k}\sigma^{-1}\imath\left(\frac{1}{k!}
\sum_{\sigma\in\SS{k}}\sum_{i=1}^{k}
f_{\beta_{\sigma(1)}}\cdots f_{\beta_{\sigma(i-1)}}
[e_{\alpha},f_{\beta_{\sigma(i)}}]
f_{\beta_{\sigma(i+1)}}\cdots f_{\beta_{\sigma(k)}}
\medspace v_{t^{2}\lambda}\right)
\end{multline*}

The term corresponding to a fixed $\sigma\in\SS{k}$ and $i=1\ldots k$
clearly vanishes unless $\alpha-\beta_{\sigma(i)}$ is a root or zero.
If $\alpha-\beta_{\sigma(i)}$ is a negative root, the corresponding
term in $\Sn_{-}$ is of degree $\leq k$ and its total contribution
an $O(1)$. On the other hand, the total contribution of the terms
for which
$$\sigma(i)\in I_{\alpha}=\{j=1\ldots k|\beta_{j}=\alpha\}$$
is
\begin{equation*}
\begin{split}
&t^{d-k}\sigma^{-1}\\
&
\left(\frac{1}{k!}
\sum_{i,\sigma:\sigma(i)\in I_{\alpha}}
\<t^{2}\lambda-\beta_{\sigma(i+1)}\cdots-\beta_{\sigma(k)},
\alpha^{\vee}\>
f_{\beta_{\sigma(1)}}\cdots f_{\beta_{\sigma(i-1)}}
f_{\beta_{\sigma(i+1)}}\cdots f_{\beta_{\sigma(k)}}\right)\\
&=
t^{d+2-k}\<\lambda,\alpha^{\vee}\>
\sigma^{-1}\left(\frac{|I_{\alpha}|}{k}\sum_{i=1}^{k}
\sigma(\frac{x_{\beta_{1}}\cdots x_{\beta_{k}}}{x_{\alpha}})\right)
+t^{-1}O(1)\\
&=
t\<\lambda,\alpha^{\vee}\>
\partial_{\alpha}
x_{\beta_{1}}\cdots x_{\beta_{k}}
+t^{-1}O(1)
\end{split}
\end{equation*}

Finally, if $\alpha-\beta_{\sigma(i)}$ is a positive root, a repetition
of the above argument shows that the net contribution is an O(1).
\eqref{eq:f} We have,
\begin{equation*}
\begin{split}
t^{d}f_{\alpha}t^{-d}
\medspace
x_{\beta_{1}}\cdots x_{\beta_{k}}
&=
t^{d-k}\sigma^{-1}(f_{\alpha}\sigma(x_{\beta_{1}}\cdots x_{\beta_{k}}))\\
&=
t^{d-k}\sigma^{-1}(\sigma(x_{\alpha}x_{\beta_{1}}\cdots x_{\beta_{k}})+r)\\
&=
tx_{\alpha}x_{\beta_{1}}\cdots x_{\beta_{k}}+O(1)
\end{split}
\end{equation*}
for some $r\in\Un_{-}$ of degree $\leq k$, where we used $f_{\alpha}
=\sigma(x_{\alpha})$ and
$$\sigma(p\cdot q)=\sigma(p)\cdot\sigma(q)+r'$$
for any $p,q\in \Sn_{-}$ where the remainder $r'\in\Un_{-}$ is
of degree $\leq\deg(p)+\deg(q)-1$.
\eqref{eq:h} follows from the fact that the eigenvalues of $h_
{\alpha}$ on $M_{t^{2}\lambda}$ lie in $t^{2}\<\lambda,\alpha^
{\vee}\>+\IZ$ and that $h_{\alpha}$ commutes with $d$ \halmos\\

Since the action of the Casimirs $C_{\alpha}$ on $M_{\mu}[\mu-\beta]$
depends polynomially on $\mu$ and irreducibility is an open condition,
it suffices to show that the set of $\mu$ for which the $C_{\alpha}$ act
irreducibly on $M_{\mu}[\mu-\beta]$ is non--empty. Let $\lambda\in\h^
{*}$ be a regular weight and choose $\mu$ of the form $t^{2}\lambda$,
where $t\in\IC^{*}$. Using the notation of lemma \ref{le:asymptotic}, it is
sufficient to show that the operators
$$t^{-2}\Ad(t^{d})f_{\alpha}e_{\alpha}
\qquad\text{and}\qquad
t^{-3}\Ad(t^{d})[f_{\alpha}e_{\alpha},f_{\beta}e_{\beta}]$$

act irreducibly on any subspace of $\Sn_{-}$ of fixed weight. Since
this is again an open condition in $t$ and, as will be shown below,
the operators at hand have a finite limit as $t\rightarrow\infty$,
it suffices to prove this for $t=\infty$. By lemma \ref{le:asymptotic},
$$\lim_{t\rightarrow\infty}
t^{-2}\Ad(t^{d})f_{\alpha}e_{\alpha}=
\<\lambda,\alpha^{\vee}\>x_{\alpha}\partial_{\alpha}$$

Let now $\alpha\neq\beta$ be positive roots. Denoting by $R\subset
\h^{*}$ the root system of $\g$, we have
\begin{equation*}
\begin{split}
[f_{\alpha}e_{\alpha},f_{\beta}e_{\beta}]
&=
f_{\alpha}e_{\alpha}f_{\beta}e_{\beta}-
f_{\beta}e_{\beta}f_{\alpha}e_{\alpha}\\
&=
f_{\alpha}f_{\beta}e_{\alpha}e_{\beta}-
f_{\beta}f_{\alpha}e_{\beta}e_{\alpha}\\
&+
\delta_{\alpha-\beta\in R}\left(
c_{\alpha,\beta}f_{\alpha}\varepsilon_{\alpha-\beta}e_{\beta}-
c_{\beta,\alpha}f_{\beta}\varepsilon_{\beta-\alpha}e_{\alpha}\right)
\end{split}
\end{equation*}

where $\varepsilon_{\gamma}=e_{\gamma}$ or $f_{\gamma}$ according
to whether the root $\gamma$ is positive or negative and the $c_{
\cdot,\cdot}$ are non--zero constants. Since
$$f_{\alpha}f_{\beta}e_{\alpha}e_{\beta}-
f_{\beta}f_{\alpha}e_{\beta}e_{\alpha}=
[f_{\alpha},f_{\beta}]e_{\alpha}e_{\beta}+
f_{\beta}f_{\alpha}[e_{\alpha},e_{\beta}]$$

we find that
\begin{equation*}
\begin{split}
[f_{\alpha}e_{\alpha},f_{\beta}e_{\beta}]
&=
\delta_{\alpha+\beta\in R}\left(
c'_{\alpha,\beta}f_{\alpha+\beta}e_{\alpha}e_{\beta}+
c'_{\beta,\alpha}f_{\beta}f_{\alpha}e_{\alpha+\beta}\right)\\
&+
\delta_{\alpha-\beta\in R}\left(
c_{\alpha,\beta}f_{\alpha}\varepsilon_{\alpha-\beta}e_{\beta}-
c_{\beta,\alpha}f_{\beta}\varepsilon_{\beta-\alpha}e_{\alpha}\right)
\end{split}
\end{equation*}

for some non--zero constants $c'_{\cdot,\cdot}$. It therefore follows
from lemma \ref{le:asymptotic} that
\begin{equation*}
\begin{split}
\lim_{t\rightarrow\infty}t^{-3}\Ad(t^{d})
[f_{\alpha}e_{\alpha},f_{\beta}e_{\beta}]
&=
\delta_{\alpha+\beta\in R}\left(
\wt{c}'_{\alpha,\beta}x_{\alpha+\beta}\partial_{\alpha}\partial_{\beta}+
\wt{c}'_{\beta,\alpha}x_{\beta}x_{\alpha}\partial_{\alpha+\beta}\right)\\
&+
\delta_{\alpha-\beta\in R}\left(
\wt{c}_{\alpha,\beta}x_{\alpha}\ol{\varepsilon}_{\alpha-\beta}\partial_{\beta}-
\wt{c}_{\beta,\alpha}x_{\beta}\ol{\varepsilon}_{\beta-\alpha}\partial_{\alpha}\right)
\end{split}
\end{equation*}

where $\ol{\varepsilon}_{\gamma}$ is now the operator $\partial
_{\gamma}$ or $x_{\gamma}$ according to whether $\gamma$ is positive
or negative and the $\wt{c}_{\cdot,\cdot},\wt{c}'_{\cdot,\cdot}$
are non--zero constants. Since the summands in the above expression
have distinct homogeneity degrees with respect to the commuting
Euler operators $x_{\alpha}\partial_{\alpha}$, it is sufficient
to show that the weight spaces of $\Sn_{-}$, \ie the subspaces
spanned by the monomials $\prod_{\alpha\succ 0}x_{\alpha}^{m_
{\alpha}}$ with $\sum_{\alpha\succ 0}m_{\alpha}\alpha$ fixed,
are irreducible under the operators
$$
x_{\alpha}\partial_{\alpha},\quad
x_{\alpha+\beta}\partial_{\alpha}\partial_{\beta},\quad
x_{\beta}x_{\alpha}\partial_{\alpha+\beta},\quad
x_{\alpha}\ol{\varepsilon}_{\alpha-\beta}\partial_{\beta},\quad
x_{\beta}\ol{\varepsilon}_{\beta-\alpha}\partial_{\alpha}$$

which is a simple enough exercise \halmos

\subsection{\sc Proof of theorem \ref{th:centraliser}}\label{ss:centraliser}
Assume that $x\in\Ugh$ commutes with the Casimirs $\calpha$. By
theorem \ref{th:Zariski}, for any $\beta\in\bigoplus_i\IN\cdot\alpha_i$,
$x$ acts as multiplication by a scalar $x(\mu,\beta)$ on the weight
space $M_{\mu}[\mu-\beta]$ of the Verma module with highest
weight $\mu\in\h^*$. Clearly, $x(\mu,\beta)$ is a polynomial in
$\mu$.

\begin{proposition}\label{pr:poly}
$x(\mu,\beta)$ is a polynomial in $\beta$.
\end{proposition}
\proof\footnote{The second author is grateful to F. Knop for outlining
the proof of proposition \ref{pr:poly}.} Note first that $x(\mu,\beta)$
is also the scalar by which $x$ acts on the weight space $M_\mu^*
[\mu-\beta]$ of the contragredient Verma module $M_\mu^*$ with
highest weight $\mu$, since the latter is isomorphic to $M_{\mu}$
for generic $\mu$.  Let now $G$ be the connected and simply--connected
complex Lie group with Lie algebra $\g$ and $N_\pm,B_\pm\subset
G$ the unipotent and solvable subgroups corresponding to $\n
_\pm$ and $\b_\pm=\n_\pm\oplus\h$ respectively. Using the action
of $G$ on the flag manifold $G/B_{-}$, identify $M_\mu^*$ as $
\b_+$--module with $\IC[N_+]\otimes \IC_\mu$ where $N_+$ is
viewed as the big cell in $G/B_-$ and $\IC_\mu$ is the
$\b_+$--module given by letting $\n_+,\h$ act by zero and $
\mu$ respectively (see, \eg \cite[\S 10.2]{FB}). This identification
yields an action of $U\g$ on $\IC[N_+]$ by differential operators
with polynomial coefficients in the generators $y_\gamma$, with
$\gamma$ a positive root, defined by $e_\alpha\cdot y_\gamma=
\delta_{\alpha,\gamma}$. The corresponding subspace of weight
$-\beta$ is spanned by the monomials $y_{\gamma_1}^{m_1}\cdots
y_{\gamma_N}^{m_N}$ such that
$m_1\gamma_1+\cdots+m_N\gamma_N=\beta$. We shall need
the following simple

\begin{lemma} If $D$ is a differential operator with polynomial
coefficients in the $y_\gamma$, the diagonal matrix entries
of $D$ in the monomial basis $y_{\gamma_1}^{n_1}\cdots
y_{\gamma_N}^{m_N}$ are polynomials in $m_1,\ldots,m_N$.
\end{lemma}
\proof It suffices to prove this for $D$ of the form $D=\partial_1^{p_1}y_1^{q_1}\cdots\partial_N^{p_N}y_N^{q_N}$.
The diagonal entries of $D$ are then clearly zero unless
$p_i=q_i$ for all $i=1\ldots N$. When that is the case
$$D y_{\gamma_1}^{m_1}\cdots y_{\gamma_N}^{m_N}=
(m_1+p_1)\cdots(m_1+1)\cdots(m_N+p_N)\cdots(m_N+1)
y_{\gamma_1}^{m_1}\cdots y_{\gamma_N}^{m_N}
$$
which is clearly a polynomial in $m_1,\ldots,m_N$ \halmos\\

In particular, if $\beta=\sum_i k_i\alpha_i$, with $k_i\in\IN$, then
$x(\mu,\beta)$ is the matrix coefficient of $x$ corresponding to
the monomial $y_{\alpha_1}^{k_1}\cdots y_{\alpha_r}^{k_r}$ and
is therefore a poynomial in $k_1,\ldots,k_r$ as claimed \halmos\\

Define now $\wt{x}\in S(\h\oplus\h)$ by $\wt{x}(\mu,\lambda)=x
(\mu,\mu-\lambda)$ so that $x$ acts as multiplication by $\wt{x}
(\mu,\lambda)$ on $M_\mu[\lambda]$. The embedding $M_{\mu
-(\mu(h_i)+1)\cdot\alpha_i}\hookrightarrow M_\mu$ valid whenever
$\mu(h_i)\in\IN$ shows that $\wt{x}$ is invariant under the $\rho
$--shifted action action of the Weyl group $W$ on the first variable
given by 
$$w\bullet \mu=w(\mu+\rho)-\rho$$
Write $\wt{x}=\sum_j z_j\cdot p_j$ with $z_j,p_j\in S\h$ and $z_j$
invariant under $W$. Regard $p_j$ as an element of $U\h$ and
let $Z_j\in Z(\Ug)$ be such that $Z_j$ acts on $M_\mu$ as multiplication
by $z_j$. Then, $x-\sum_j Z_j\cdot p_j$ acts as zero on all Verma
modules and is therefore equal to zero since these separate elements
in $\Ug$ \halmos

\end{document}